\documentclass[a4paper,12pts,oneside]{amsart}
\usepackage[top=1in,bottom=1.1in,left=1.35in,right=.95in]{geometry}
\usepackage{amssymb}
\usepackage[utf8]{inputenc}
\usepackage[T1]{fontenc}
\usepackage[french,english]{babel}
\usepackage[all]{xy}
\usepackage{soul}
\usepackage{ulem}
\usepackage{shuffle}
\usepackage{hyperref}
\usepackage{stmaryrd}
\usepackage{graphicx}
\usepackage{array}
\usepackage{xcolor}
\usepackage{slashed}
\usepackage{amsthm}
\usepackage{amsmath}
\usepackage{amssymb}
\usepackage{mathrsfs}
\newtheorem{Lemme}{Lemme}[section]
\newtheorem{Prop}{Proposition}[section]
\newtheorem{Def}{Definition}[section]
\newtheorem{Rem}{Remark}[section]
\newtheorem{Thm}{Theorem}[section]
\newtheorem{Cor}{Corollary}[section]
\newtheorem*{Not}{Notation}

\newtheorem{Thmintro}{Theorem}

\newtheorem{Example}{Example}[section]

\newcommand{\Preuve}{\noindent\textbf{Proof:\ }}
\newcommand{\eb}{\vspace{-1em}\begin{flushright}$\Box$\end{flushright}}
\newcommand{\eq}[1][r]
   {\ar@<-3pt>@{-}[#1]
    \ar@<-1pt>@{}[#1]|<{}="gauche"
    \ar@<+0pt>@{}[#1]|-{}="milieu"
    \ar@<+1pt>@{}[#1]|>{}="droite"
    \ar@/^2pt/@{-}"gauche";"milieu"
    \ar@/_2pt/@{-}"milieu";"droite"}

\newcommand{\cF}{{\mathcal{F}}}
\newcommand{\tcF}{{\widetilde{\mathcal{F}}}}
\newcommand{\cT}{{\mathcal{T}}}
\newcommand{\cU}{{\mathcal{U}}}
\newcommand{\cV}{{\mathcal{V}}}

\newcommand{\cH}{{\mathcal{H}}}

\newcommand{\tM}{{\widetilde{M}}}

\newcommand{\Hol}{{\mathrm{Hol}}}

\newcommand{\ta}{{\tilde{\frak{a}}}}

\newcommand{\ha}{\hat{\frak{a}}}
\newcommand{\fa}{\frak{a}}
\newcommand{\hM}{{\widehat{M}}}
\newcommand{\hcT}{{\widehat{\mathcal{T}}}}
\newcommand{\hcF}{{\widehat{\mathcal{F}}}}

\begin{document}
\selectlanguage{english}
\title{A generalization of Molino's theory and equivariant basic $\hat{A}$-genus characters}
\author{Wenran Liu}
\address{Wenran Liu, SIAE Civil Aviation University of China, 2898 Jinbei Street Dongli District, 300300 Tianjin, China }
\email{wrliu@cauc.edu.cn}
\thanks{This research project is supported by supported by the Fundamental Research Funds for the Central Universities, SCUT}
\keywords{Riemannian foliations, Molino's theory, Killing foliations, equivariant cohomology, equivariant $\hat{A}$-genus character}
\date{}

\begin{abstract}
Molino's theory is a mathematical tool for studying Riemannian foliations. In this paper, we propose a generalization of Molino's theory with two Riemannian foliations. For this purpose, the projection of foliation with respect to a fibration is discussed. The generalization results in an equivariant basic cohomological isomorphism in case of Killing foliation. It is a generalization of results given by Goertsches and Töben. We also give a geometric realization of the cohomological isomorphism through equivariant basic $\hat{A}$-genus characters, who play a prominent role in calculating the index of an elliptic operator by Atiyah-Singer's index formula. 
\end{abstract}

\maketitle

\tableofcontents
\normalsize
\section{Introduction}
The notion of bundle-like metric on foliated manifold was first introduced by Reinhart \cite{R}. A foliated manifold equipped with a bundle-like metric is called Riemannian foliation. Molino's theory \cite{Mol} is a mathematical tool for studying Riemannian foliations. Roughly, to each transversely oriented Riemannian foliation $(M,\cF)$ of codimension $q$, Molino associated an oriented manifold $W$ equipped with an action of the orthogonal group $SO(q)$. Let $E=E^+\oplus E^-$ be a foliated vector bundle over $M$ and $\mathcal{D}_b^{E+}:C^\infty_{\cF\text{-bas}}(M,E^+)\rightarrow C^\infty_{\cF\text{-bas}}(M,E^-)$ be a transversely elliptic basic differential operator on basic sections. Consider its basic index $\mathrm{Index}_b(\mathcal{D}_b^{E+})\in\mathbb{Z}$. It is an open question to obtain a cohomological formula as Atiyah-Singer index theorem for this invariant. Some partial responses have been proposed in the literature \cite{BKR,GL,Ka,PV}.\\

In the 1990s, El Kacimi-Alaoui proposed to tackle this problem by using Molino's theory. The main result of \cite{Ka} is the association of a $SO(q)$-equivariant operator $\mathcal{D}^{\mathcal{E}+}:C^\infty(W,\mathcal{E}^+)\rightarrow C^\infty(W,\mathcal{E}^-)$ to $\mathcal{D}_b^{E+}$ and the relation $\mathrm{Index}_b(\mathcal{D}_b^{E+})=\dim\big[\mathrm{Index}_{SO(q)}(\mathcal{D}^{\mathcal{E}+})\big]^{SO(q)}$ where $\big[-\big]^{SO(q)}$ denotes the $SO(q)$-invariant part of the equivariant index of $\mathcal{D}^{\mathcal{E}+}$. Following this valuable identity, we investigate the acquisition of $\mathrm{Index}_b(\mathcal{D}_b^{E+})$ from the information of $\mathcal{D}^{\mathcal{E}+}$. \\

Let $H^\bullet(M,\cF)$ be the basic de Rham cohomology of the foliated manifold $(M,\cF)$. If $(M,\cF)$ is a Killing foliation (a priori Riemannian), there is an Abelian algebra $\fa$, acting transversely on $(M,\cF)$. The orbit of the leaves of $\cF$ under the action of $\fa$ is exactly the leaf closure. Let $H_\fa^\bullet(M,\cF)$ be the $\fa$-equivariant basic de Rham cohomology with polynomial coefficients, see \cite{GS,KT}. Goertsches and Töben \cite{GT} proved a cohomological isomorphism $H_\fa^\bullet(M,\cF)\simeq H_{SO(q)}^\bullet(W)$. Duflo and Vergne \cite{DKV} extended the equivariant cohomology to the case of $C^\infty$-coefficients. \cite{Liu} showed the geometric realization of the cohomological isomorphism through the equivariant basic Chern characters. This is the first step following Kacimi's proposal. \\

Our next step is to establish a geometric realization through the equivariant basic $\hat{A}$-genus characters who play an important role in Atiyah-Singer index formula. In this work, we construct it in a larger framework, that is, proceed to generalize Molino's theory for foliated manifold $(M,\cF)$ equipped with a second foliation $\cF^\prime$ with $\cF\subset\cF^\prime$ (each leaf of $\cF$ lies in a leaf of $\cF^\prime$). A natural question is: can $\cF^\prime$ be associated to a foliation $\cF_W$ on $W$. It is an affirmative answer in case that $\cF^\prime$ is a Riemannian foliation. A counterexample is presented in Remark \ref{Rem8}, thus the Riemannian assumption of $\cF^\prime$ is essential. \\

\noindent Followed by the generalization of Molino's theory, we obtain a cohomological isomorphism
$$
H^\infty_{\fa_\perp}(M,\cF)\simeq H^\infty_{so(\hat{q})}(W,\cF_W),
$$
where $\fa_\perp$ is a Lie subalgebra of $\fa$ which acts transversely on $\cF^\prime$ and $\hat{q}$ a constant which relates to the codimensions of $\cF$ and $\cF^\prime$. Notice that if $\cF^\prime=\cF$, $\cF_W$ degenerates to points and $\fa_\perp=\fa$, that is \cite{GT} Proposition 4.9. Similar to the definition of equivariant basic Chern characters, the equivariant basic $\hat{A}$-genus character of normal bundle $\hat{A}_{\fa_\perp}(\nu\cF)$ \big(resp. $\hat{A}_{so(\hat{q})}(\nu\cF_W)$\big) is well-defined. Finally, we establish a geometric realization through equivariant basic $\hat{A}$-genus characters of normal bundle.\\

\noindent The main result of this article is as follows.
\begin{Thmintro}\label{Thmintro1}
Let $\cF$ be the simple foliation of a fibration $M\rightarrow W$ of two compact manifolds. If $\cF^\prime$ is a foliation on $M$ with $\cF\subset \cF^\prime$, then the foliation $\cF^\prime$ projects to a foliation $\cF_W$ on $W$. Furthermore, if $\cF^\prime$ is a Riemannian foliation, $\cF_W$ also is.
\end{Thmintro}
\begin{Thmintro}\label{Thmintro2}
Let $M$ be a foliated manifold equipped with Riemannian foliations
$\cF, \cF^\prime$ satisfying $\cF\subset \cF^\prime$. Let $\widehat{M}\!=\!O(\cF^\prime\slash\cF)$ be the orthonormal frame bundle of $\cF^\prime\slash\cF$ and $\hcF^\prime$ be the pullback foliation on $\hM$. 
\begin{enumerate}
\item $\widehat{M}$ is a principal $\cF$-bundle with a lifted Riemannian foliation $\hcF$. The foliation $\hcF$ is partial transversely parallelizable to $\hcF^\prime$;
\item $\tM\!=\!O(\nu\widehat{\cF})$ is a principal $\cF$-bundle with lifted Riemannian foliations $\tcF,\tcF^\prime$. The foliations $\tcF,\tcF^\prime$ are transversely parallelizable in $\tM$ and $\tcF$ is partial transversely parallelizable to $\tcF^\prime$;
\item $\tcF^\prime$ projects to a Riemannian foliation $\cF_W$ on $W$. 
\end{enumerate}
See the diagram below.
\begin{equation*}
\xymatrix{
&(\tM,\tcF,\tcF^\prime)\ar[r]^-{\pi}\ar[d]^-{\widetilde{p}}&(W,\cF_W)\\
(M,\cF,\cF^\prime)&(\hM,\hcF,\hcF^\prime)\ar[l]^-{\widehat{p}}&
}
\end{equation*}  
\end{Thmintro}
\begin{Thmintro}\label{Thmintro3}
Let $M$ be a foliated manifold equipped with a Killing foliation $\cF$ with Molino algebra $\fa$ and a Riemannian foliation $\cF^\prime$ satisfying $\cF\subset \cF^\prime$. Let $\fa_\perp$ be the Lie subalgebra orthogonal to $\cF^\prime$ of $\fa$. Then, there is a cohomological isomorphism as follows. Furthermore, the equivariant basic $\hat{A}$-genus characters of normal bundle give a geometric realization in following way
\[
\begin{array}{ccc}
H^\infty_{\fa_\perp}(M,\cF^\prime)&\simeq&H^\infty_{so(\hat{q})}(W,\cF_W)\\ \\[-0.5em]
\hat{A}_{\fa_\perp}(\nu\cF^\prime)&\simeq&\hat{A}^{-1}_{so(\hat{q})}\cdot\hat{A}_{so(q)}(\nu\cF_W),
\end{array}
\]
where \small $\hat{A}_{so(\hat{q})}=\det{}^{\!\frac{1}{2}}\left[\dfrac{(\mathrm{ad}X)/2}{\sinh((\mathrm{ad}X)/2)}\right]$\normalsize  is the function of $so(\hat{q})$ with a constant $\hat{q}$.
\end{Thmintro}

\noindent In Section 2, we review some basic definitions and properties of foliation, foliated bundle, groupoid action, Assumption, etc.\medskip\ \\
\noindent In Section 3, we propose some properties of manifold with two foliations. In this section, some properties are regardless of the Riemannian assumption of foliation.\medskip\ \\
\noindent In Section 4, we discuss the projection of foliation, it is a reversed work of \cite{Mac} Section 1.6. This part will be applied in Section 5 for associating a foliation on the basic manifold $W$.\medskip\ \\
\noindent In Section 5, we generalize Molino's theory. In order to lift the second foliation $\cF^\prime$ and keep compatible with $\cF$, two lifts are necessary (only once as Molino's is insufficient). As showed in the diagram of Theorem \ref{Thmintro2}, two lifts result respectively in $(\hM,\hcF,\hcF^\prime)$ and $(\tM,\tcF,\tcF^\prime)$. When $\cF^\prime=\cF$, the first lift is degenerated, back to Molino's theory.\medskip\ \\
\noindent Section 6 is devoted to the generalization of cohomological isomorphism proposed in \cite{GT} Proposition 4.9. In this section, the foliation $\cF$ (resp. $\cF^\prime$) is supposed to be Killing (resp. Riemannian). Remark that only the subalgebra orthogonal to $\cF^\prime$ of Molino algebra acts transversely to $\cF^\prime$.\medskip\ \\
\noindent In Section 7, the equivariant basic $\hat{A}$-genus characters are considered to give a geometric realization of the cohomological isomorphism obtained in Section 6. \cite{Liu} Theorem 5.4.1 also \cite{Liu1} Theorem 6 showed the geometric realization by the equivariant basic Chern characters. Notice that equivariant basic Chern characters are perfectly identified. A slight difference for $\hat{A}$-genus character is that their identification is up to a function $\hat{A}_{so(\hat{q})}$ of $so(\hat{q})$. \\

\noindent {\bf\em Acknowledgements:} The author is grateful to Moulay Tahar Benameur and Alexandre Baldare for various discussions on this project. The author is also grateful to the support of the Fundamental Research Funds for the Central Universities, SCUT. Project No. 3122020074. \\
\noindent The author thankfully welcome and accept remarks, comments, doubts and critiques about this paper. \\

\noindent Contact: \textbf{wrliu@cauc.edu.cn} also \textbf{wenran.liu@gmail.com}

\section{Preliminaries}
\noindent The needed basic concepts from foliation theory can be found in \cite{CC,HH1,HH2}.
\subsection{Basic definitions}
Let $\cF$ be a smooth foliation of codimension $q$ on a compact manifold $M$. Let $T\cF$ denote the vector subbundle of vectors tangent to the leaves and $\nu\cF=TM\slash T\cF$ its normal bundle. In this paper, assume that all foliations are transversely orientable, i.e. $\nu\cF$ is orientable. Denote by $\Gamma(\nu\cF)$ the Lie algebra of the sections of $\nu\cF$, by $\mathcal{L}$ the \textit{Lie derivative} on corresponding spaces. 
\begin{Def} A \textbf{transverse metric} on $(M,\cF)$ is a symmetric positive $C^\infty(M)$-bilinear form on $\Gamma(\nu\cF)$ such that
$$
\mathcal{L}_X g=0 
$$
for any vector field $X$ on $M$ tangent to $\cF$. $\cF$ is called a \textbf{Riemannian foliation} if there exist a transverse metric. Simply, we say that $\cF$ is Riemannian. 
\end{Def}
The word ``transverse'' means the holonomy invariance, that is, each holonomy transformation on $\nu\cF$ is an isometry. We assume that the reader is familiar with the notion of \textbf{holonomy groupoid}, see \cite{Moe} Section 5.2. Let $\Hol(M,\cF)$ be the holonomy groupoid of $(M,\cF)$. We denote by $s$ and $r$ for the source and target map on all groupoids.
\begin{Prop}[\cite{Moe} Proposition 5.6 Example 5.8(7)] The holonomy groupoid of a foliation has a natural Lie groupoid structure. Furthermore, if the foliation is Riemannian, its holonomy groupoid is Hausdorff.
\end{Prop}
For the notion of groupoid action on a fibration, see \cite{Mac} Section 1.6. When discuss a groupoid action on a fibration, in particular on a vector or principal bundle, we adopt a convention that the groupoid is \textbf{Lie} and \textbf{Hausdorff}. Here are respectively the definition of $\Hol(M,\cF)$-action on principal (resp. Euclidean) bundle.
\begin{Def}[\cite{Liu} Definition 4.1.2,4.1.3]
\ \\
\vspace{-1em}
\begin{enumerate}
\item A $\Hol(M,\cF)$-\textbf{equivariant principal bundle} $P\rightarrow M$ is a principal bundle $P$ equipped with a groupoid action of $\Hol(M,\cF)$ such that for each $\gamma\in \Hol(M,\cF)$, the map
$$
T_{\gamma}:P\vert_{s(\gamma)}\rightarrow P\vert_{r(\gamma)}
$$
is equivariant under the structure group of $P$;
\item A $\Hol(M,\cF)$-\textbf{equivariant Euclidean (resp. hermitian) bundle} $E\rightarrow M$ is a Euclidean (resp. hermitian) bundle $E$ equipped with a groupoid action of $\Hol(M,\cF)$ such that for each $\gamma\in \Hol(M,\cF)$, the map
$$
T_{\gamma}:E\vert_{s(\gamma)}\rightarrow E\vert_{r(\gamma)}
$$
is an isometry. 
\end{enumerate}
For simplicity, we call $P$ (resp. $E$) is $\Hol(M,\cF)$-equivariant.
\end{Def}
\begin{Rem}[\cite{Liu} Definition 4.1.4]\label{Rem2.1} If a Euclidean (resp. hermitian) bundle $E$ is $\Hol(M,\cF)$-equivariant, its orthonormal (resp. unitary) frame bundle $O(E)$ (resp. $U(E)$) is naturally $\Hol(M,\cF)$-equivariant.
\end{Rem}
\begin{Def}[Foliated bundle]
A \textbf{foliated principal bundle} $\pi:(P,\cF_P)\rightarrow (M,\cF)$ is a principal bundle $P$ equipped with a foliation $\cF_P$ such that 
\begin{itemize}
\item[$(1)$] $\cF_P$ is invariant under the action of the structure group of $P$;
\item[$(2)$] for each $p\in P$, the projection $T\pi\vert_p:T\cF_P\vert_p\rightarrow T\cF\vert_{\pi(p)}$ is an isomorphism. 
\end{itemize}
A \textbf{foliated Euclidean (resp. hermitian) bundle} is a Euclidean (resp. hermitian) bundle associated to a foliated principal bundle.
\end{Def}
\begin{Def}[$\cF$-bundle]
A \textbf{principal $\cF$-bundle} $(P,\cF_P)$ is a principal foliated bundle equipped with a basic connection $\omega$. That is, a connection $\omega$ on $P$ satisfies 
$$
i_X \omega=\mathcal{L}_X \omega=0
$$
for any vector fields $X$ tangent to $\cF_P$, where $i$ is the contraction. A \textbf{Euclidean} (resp. \textbf{hermitian}) \textbf{$\cF$-bundle} is a Euclidean (resp. hermitian) bundle associated to a principal $\cF$-bundle.
\end{Def}
\noindent Foliated bundles are precisely discussed in \cite{KT}. The following two propositions summarize some basic properties.
\begin{Prop}
A groupoid action of $\Hol(M,\cF)$ on a principal bundle $P$ defines a foliation $\cF_P$ on $P$ such that $(P,\cF_P)$ is a foliated principal bundle. 
\end{Prop}
\begin{Prop}[\cite{Liu} Proposition 4.1.4]\label{Prop-3}
Let $\Hol(M,\cF)\ltimes P$ be the semi-product groupoid. There is a Lie groupoid isomorphism   
$$
\Hol(M,\cF)\ltimes P\simeq \Hol(P,\cF_P).
$$
\end{Prop}
\subsection{Holonomy groupoid closure}
A \textbf{transversal} $\cT$ of $\cF$ is a submanifold of $M$ of dimension $q$ such that $\forall\ m\in\cT$, $T_m\cT\cap T\cF\vert_m=\{0\}$. A \textbf{complete transversal} is an immersed (not necessarily connected) submanifold, which intersects any leaf in at least one point. Remark that a complete transversal section always exists \cite{Moe}.\\
\indent Let $\cT$ be a complete transversal of $\cF$, we identify the normal bundle $\nu\cF$ with the tangent bundle $T\cT$ and the transverse metric $g$ with a Riemannian metric on $\cT$. Let $\Hol(M,\cF)_\cT^\cT$, $\Hol(P,\cF_P)_{P_\cT}^{P_\cT}$ be respectively the étale restricted holonomy groupoid on $\cT$ and $P\vert_\cT$ (Denote simply $P_\cT$ for $P\vert_\cT$). It follows an isomorphism from Proposition \ref{Prop-3}
$$
\Hol(M,\cF)_\cT^\cT\ltimes P_\cT\simeq \Hol(P,\cF_P)_{P_\cT}^{P_\cT}.
$$
\noindent Let
\begin{equation}\label{J1t}
J^1(\cT)=\big\{ (x,y,A)\mid x,y\in\cT,\ A:T_x\cT\rightarrow T_y\cT  \ \text{is}\ g\ \text{isometric}\big\}
\end{equation}
be the \textbf{isometric} $1$-jet groupoid of $\cT$. Endowed with \textbf{$J^1$-topology}, $J^1(\cT)$ is \textbf{Lie Hausdorff} and \textbf{proper}. Consider 
$$
j:\Hol(M,\cF)_\cT^\cT\rightarrow J^1(\cT)
$$ 
the injective groupoid homomorphism which maps a germ of local diffeomorphism to its derivative. It is injective because a germ of local isometry is completely determined by its $1$-jet. For detail, see \cite{GL,Liu}.   
\begin{Def}
The groupoid closure $\overline{\Hol(M,\cF)_\cT^\cT}$ is the closure of $j\Big(\Hol(M,\cF)_\cT^\cT\Big)$ in $J^1(\cT)$. 
\end{Def}
\noindent In \cite{Liu}, the following assumption is introduced, also Assupmtion 3 in \cite{Liu1}. \\ \\
\textbf{Assumption: The groupoid action of $\Hol(M,\cF)_\cT^\cT$ on $P_\cT$ extends to a groupoid action of $\overline{\Hol(M,\cF)_\cT^\cT}$ on $P_\cT$.}
\begin{Prop}[\cite{Liu} Theorem 4.1] \label{Prop-2}
Under \textbf{Assumption}, $(P,\cF_P)$ is a principal $\cF$-bundle.
\end{Prop}
\begin{Prop}[\cite{Liu} Appendix]\label{Prop-1.1} $\cF_P$ is a Riemannian foliation if and only if $(P,\cF_P)$ is a principal $\cF$-bundle.
\end{Prop}
\begin{Rem}\label{Liftmetric}
A basic connection $\omega$ on $P$ is equivalent to a transverse metric on $(P,\cF_P)$. Indeed, define the transverse metric $g_P=\pi^{-1}g\oplus g_{\cV}$ where $\pi^{-1}$ is the horizontal lift with respect to $\omega$ and $g_{\cV}$ is the canonical metric on vertical subbundle. 
\end{Rem}
\noindent Let $J^1(P_\cT)$ be the isometric $1$-jet groupoid of $P_\cT$ with respect to $g_P$. 

\begin{Prop}\label{Prop-1.5}
The following groupoid homomorphism is well-defined and injective
\[
\begin{array}{rcl}
k:\overline{\Hol(M,\cF)_\cT^\cT}\ltimes P_\cT&\rightarrow& J^1(P_\cT)\\
(A,p)&\mapsto&\big(p,A\!\cdot\! p,\widetilde{A}\oplus \mathrm{id}_{\cV}),
\end{array}
\]
where $\widetilde{A}$ is the horizontal lift of $A$ and $\mathrm{id}_{\cV}$ is the identity map on vertical subbundle.
\end{Prop}
\noindent See the commuting diagram
$$
\xymatrix{
\Hol(M,\cF)_\cT^\cT\ltimes P_\cT \ar[d]^{j}\ar[r]^-{\simeq} &\Hol(P,\cF_P)_{P_\cT}^{P_\cT}\ar[d]^-{j_P}\\
\overline{\Hol(M,\cF)_\cT^\cT}\ltimes P_\cT\ar[r]^-{k}&J^1(P_\cT) 
}
$$
where $j_P$ is the injective groupoid homomorphism similar to $j$. Take groupoid closure, followed by injectivity of $j$, $j_P$ and $k$, we obtain the following proposition.
\begin{Prop}\label{Prop-1}
There is a Lie groupoid isomorphism
$$
\overline{\Hol(M,\cF)_\cT^\cT}\ltimes P_\cT\simeq \overline{\Hol(P,\cF_P)_{P_\cT}^{P_\cT}}.
$$
\end{Prop}
Resulted from \cite{Sal1,Sal2}, next propositions reveal the relation between the Molino sheaf of a Riemannian foliation and the groupoid closure.  
\begin{Prop}
$\forall m\in \cT$, each element of $\overline{\Hol(M,\cF)_\cT^\cT}$ near enough to $1_m$ is written by
$(m,e^{tY_\cT}\!\cdot\! m, T_m e^{tY_\cT})$ for a local transverse field $Y$ belonging to the Molino sheaf (identified with a local field $Y_\cT$ on $\cT$). Conversely, each element in this form for $t$ small enough belongs to $\overline{\Hol(M,\cF)_\cT^\cT}$.
\end{Prop}

\begin{Prop}\label{Prop0}
The groupoid action of $\overline{\Hol(M,\cF)_\cT^\cT}$ on $P_\cT$ induces an identification of Molino sheaf on $(M,\cF)$ and $(P,\cF_P)$.
\end{Prop}

\section{Manifold with two foliations}
\noindent Let $M$ be a compact manifold equipped with two foliations $\cF,\cF^\prime$ with $\cF\subset\cF^\prime$, i.e. each leaf of $\cF$ lies in a leaf of $\cF^\prime$, equivalently $T\cF$ is a subbundle of $T\cF^\prime$. Before giving definitions and properties, let us review an example of non-Riemannian foliation presented in \cite{Liu1} Example 2.
\begin{Example}\label{Ex1}
Let $M\!=\!S^1\times S^1\times S^1$ with local coordinates $(x,y,z)$. Take $f$ a non-constant function on $S^1$. Two foliations are defined by
$\cF\!=\!\big<\partial\slash\partial x\!+\!f(y)\partial\slash\partial z \big>$ and $\cF^\prime\!=\!\big<\partial\slash\partial x,\partial\slash\partial z\big>$. 
\end{Example}
Clearly, $\cF^\prime$ is trivial foliation, thus Riemannian. On each leaf of $\cF^\prime$, $\cF$ is Kronecker foliation, thus Riemannian. But $\cF$ is a non-Riemannian foliation on $M$. From which we know that the Riemannian property of foliation is non-hereditary. Hence, supposing the Riemannian assumption of $\cF$ is essential. From now on, $\cF$ is a \textbf{Riemannian foliation}.  
\begin{Prop}\label{Lemme1}
Let $\cT$ be a transversal of $\cF$. For each $m\in \cT$, there exist a transversal $\cU$ of $\cF^\prime$ at $m$ with $\cU\subset \cT$. 
\end{Prop}
\Preuve Let $U$ be a simple open set in $M$. Consider $\cT$ to be the transversal denoted by $T_U$ in \cite{Mol} Page 21. Suppose that $\cT$ meets each plaque of $\cF$ in $U$ at a unique point. Let $\pi:U\rightarrow \cT$ be the projection along leaves of $\cF$. Take $m_0\in M$ such that $\pi(m_0)=m$. Let $\cU_0$ be a transversal of $\cF^\prime$ at $m_0$ located in $U$ such that $\cU_0$ meets each plaque of $\cF^\prime$ in $U$ at a unique point. \smallskip\ \\
Since $\cF\subset \cF^\prime$, $\cU_0$ meets also each plaque of $\cF$ in $U$ at a unique point if it does. Thus, the restricted map $\pi: \cU_0\rightarrow \pi(\cU_0)$ is a bijection. Each vector tangent to $\cU_0$ is transverse to $T\cF^\prime$, a priori transverse to $T\cF$. So $\forall\ m_0\in\cU_0$, the tangent map $T\pi\vert_{m_0}$ is an isomorphism. Hence, $\cU=\pi(\cU_0)$ is an immersed submanifold of $\cT$. In fact, $\pi:\cU_0\rightarrow \cU$ is an \'{e}tale map. \smallskip \ \\ 
If $x\in\cU$, $v\in T\cF^\prime\vert_x$, there exist a curve $c:(-\epsilon,\epsilon)\rightarrow N$ such that $\dot{c}(0)=v$ where $N$ is the leaf of $\cF^\prime$ passing by $x$. $\pi$ defines a lift $\tilde{c}$ of $c$ in $\cU_0$ with $\tilde{c}(0)=x_0$. As for each $t$, $c(t)$ and $\tilde{c}(t)$ lie in same leaf of $\cF$, also in same leaf of $\cF^\prime$, the curve $\tilde{c}$ lies in $N$. Hence, the lifted vector $\tilde{v}=\dot{\tilde{c}}(0)\in T\cF^\prime\vert_{x_0}$. As $T_{x_0}\cU_0\cap  T\cF^\prime\vert_{x_0}=\{0\}$ induces $T_x\cU\cap T\cF^\prime\vert_x=\{0\}$, $\cU$ is a transversal of $\cF^\prime$ with $\cU\subset \cT$. \eb
\begin{Def}
We call $\cU$ a \textbf{subtransversal} of $\cT$ at $m$.  
\end{Def}
Followed by the proof of Proposition \ref{Lemme1}, we have the following corollary.
\begin{Cor}\label{subtransversal}
For each complete transversal $\cT$ of $\cF$, there exist a complete transversal $\cU$ of $\cF^\prime$ such that $\cU\subset \cT$. We call $\cU$ a complete subtransversal of $\cT$.
\end{Cor}
\begin{Cor}\label{Lemme2}
To each $\gamma\in \Hol(M,\cF)$, any representative $f:(U,x)\rightarrow (V,x)$ where $U,V$ are transversals at source $x$ and target $y$. The image of each subtransversal of $U$ at $x$ by $f$ is a subtransversal of $V$ at $y$.
\end{Cor}
\smallskip\ 
The normal bundle $\nu\cF$ is $\Hol(M,\cF)$-equivariant. Denote by $\xi\!=\!\cF^\prime\slash\cF$, equipped with the restricted transverse metric of $\nu\cF$. From the proof of Proposition \ref{Lemme1}, we know that $\xi$ is stable under the action of $\Hol(M,\cF)$. 
\begin{Prop}\label{Prop1}\ \\
\vspace{-1em}
\begin{itemize}
\item[$(i)$] The Euclidean bundle $\xi$ is $\Hol(M,\cF)$-equivariant. 
\item[$(ii)$] There is a natural holonomy groupoid homomorphism $\Hol(M,\cF)\rightarrow \Hol(M,\cF^\prime),\ \gamma\mapsto \underline{\gamma}$. 
\end{itemize}
\end{Prop}
\noindent Each $\cF$-transverse metric $g$ on $\nu\cF$ defines an orthogonal decomposition 
\begin{equation}\label{eq:ortho}
\nu\cF=\xi\oplus \xi_\perp. 
\end{equation}
A natural isomorphism is $\nu\cF^\prime\simeq \xi_\perp$. The restriction of $g$ on $\xi_\perp$ is identified with a $\cF$-transverse metric on $\nu\cF^\prime$, which summarizes to the following corollary.
\begin{Cor}\label{Cor1}
The normal bundle $\nu\cF^\prime$ is $\Hol(M,\cF)$-equivariant. 
\begin{itemize}
\item[$(i)$] For each $\gamma\in \Hol(M,\cF)$ and $v\in \nu\cF^\prime\vert_{s(\gamma)}$, we have
\begin{equation}\label{eq:02}
\underline{\gamma}\cdot v=(\gamma\cdot \tilde{v})_\perp,
\end{equation}
where $\tilde{v}\in\xi\vert_{s(\gamma)}$ is an arbitrary representative of $v$, where $\perp:\nu\cF\rightarrow \nu\cF^\prime$ is the projection;
\item[$(ii)$] For each $\gamma\in \Hol(M,\cF)$ and $u\in \nu\cF\vert_{s(\gamma)}$, we have
$$
\underline{\gamma}\cdot u_\perp=(\gamma\cdot u)_\perp,
$$
\end{itemize}
\end{Cor}
The previous properties hold without Riemannian assumption on $\cF^\prime$. Next, suppose that $\cF,\cF^\prime$ are both \textbf{Riemannian} foliations equipped with transverse metrics $g_1$ and $g_2$, respectively.
\begin{Def}
Followed by relation (\ref{eq:ortho}), we define the transverse metric $g_1\vee g_2$ on $\nu\cF$ by
\begin{equation}\label{eq:combine}
g_1\vee g_2={g_1}\vert_{\xi}\oplus g_2.
\end{equation}
\end{Def}
\noindent The following propositions result from the orthogonality of (\ref{eq:ortho}) under the action of $\Hol(M,\cF)$.
\begin{Prop}\label{Prop:metric}
The metric $g_1\vee g_2$ on $\nu\cF$ is $\cF$-transverse.
\end{Prop}
\begin{Prop}\label{Prop1.1}
The orthonormal frame bundle $O(\nu\cF^\prime)$ (resp. $O(\xi)$) with respect to $g_2$ (resp. $g_1\vert\xi$) on $\nu\cF^\prime$ (resp. $\xi$) satisfies \textbf{Assumption} for the action of $\Hol(M,\cF)$ in relation (\ref{eq:02}).
\end{Prop}
\section{Projection of Foliation}
In this section, let $\cF$ be the \textbf{simple foliation} of a fibration $\pi:M\rightarrow W$ of two compact manifolds. Let $\cF^\prime$ be a  foliation on $M$ (\textbf{not necessarily Riemannian}) with $\cF\subset \cF^\prime$.  In this section, we prove that the foliation $\cF^\prime$ ``projects'' to a foliation $\cF_W$ on $W$. This is an reversed work of the lift of foliation in \cite{Mac} Section 1.6. The main line of proof is to construct a Lie groupoid with base $W$ whose orbit naturally defines $\cF_W$. \medskip \ \\
We know that the holonomy groupoid $\Hol(M,\cF)$ of simple foliation is isomorphic to the pair of points on same fiber, that is
$$
\Hol(M,\cF)\simeq M\times_\pi M=\big\{(x,y)\mid \pi(x)=\pi(y)\big\}.
$$ 
Let $1_x^y$ denote simply the element $(x,y)$. Notice that $\Hol(M,\cF)\rightarrow \Hol(M,\cF^\prime)$ is injective.
Denote $\underline{1_x^y}\!\in\! \Hol(M,\cF^\prime)$ the image.
\begin{Def}
Let $\cdot$ denote the product on $\Hol(M,\cF^\prime)$. On the arrow set of $\Hol(M,\cF^\prime)^{(1)}$, we define a product $*$: for $\gamma,\delta\in\Hol(M,\cF^\prime)$ with $\pi(r(\delta))=\pi(s(\gamma))$, 
$$
\gamma * \delta=\gamma \cdot \underline{1_{r(\delta)}^{s(\gamma)}}\cdot \delta.
$$
\end{Def}
\begin{Def}
Let $\gamma,\delta\in\Hol(M,\cF^\prime)$. $\gamma$ is equivalent to $\delta$, denoted by $\gamma\sim \delta$, if
\begin{enumerate}
\item their sources and targets are on the same fiber of $\pi$,
\item  $\delta*\gamma^{-1}=\underline{1_{r(\gamma)}^{r(\delta)}}$
\end{enumerate} 
\end{Def}
\noindent It is easy to check that $*$ is well-defined and $\sim$ is an equivalence relation. Let $[\gamma]$ denote the equivalence class of $\gamma$.
\begin{Def}[Holonomy groupoid projection]
A groupoid $\cH$ is defined by following data
\begin{enumerate}
\item the base manifold $\cH^{(0)}=W$;
\item the arrow set $\cH^{(1)}=\Hol(M,\cF^\prime)^{(1)}\slash_\sim$, that is the quotient set under $\sim$;
\item the source and target map $s([\gamma])=\pi(s(\gamma))$,$r([\gamma])=\pi(r(\gamma))$;
\item the product $[\gamma][\delta]=[\gamma*\delta]$ for $[\gamma],[\delta]\in \cH^{(1)}$ with $s([\gamma])=r([\delta])$;
\item the unit inclusion $i:W\rightarrow \cH^{(1)}, w\mapsto [\underline{\mathrm{1}^x_x}]$ with $x\in\pi^{-1}(w)$.
\end{enumerate}
\end{Def}
\begin{Prop}\label{Prop4}
The groupoid $\cH$ has a Lie groupoid structure.
\end{Prop}
\Preuve In this proof, let $\dim\cF^\prime\!=\!p$, $\mathrm{Codim}\cF^\prime\!=\!q$. The proof contains four parts.\\
\underline{Part I: Local charts for $\cH$} To each $[\gamma]\in \cH$ with $\gamma\in \Hol(M,\cF^\prime)$ be a representative. By \cite{Moe} Proposition 5.6, local foliated charts at $x\!=\!s(\gamma)$ and $y\!=\!r(\gamma)$ are written by
$$
\varphi\!:\!U\rightarrow A\,\!\times \!C\subset \mathbb{R}^{p}\times \mathbb{R}^{q},\ \  \psi\!:\!V\rightarrow B\,\!\times\,\! D \subset \mathbb{R}^{p}\times \mathbb{R}^{q},
$$ 
with $\varphi(x)=(a,c)$, $\psi(y)=(b,d)$. A plaque of $\cF^\prime$ at $x$ is $\varphi^{-1}(A\times\{c\})$. A local chart of $\gamma$ is given by $f:A\times B\times C\rightarrow \Hol(M,\cF^\prime)^{(1)}$. \smallskip\ \\
As $\cF\subset\cF^\prime$, $\varphi^{-1}(A\times\{c\})$ is equipped with a foliated chart of $\cF$, then there exist a small open $A^\prime$ at $a$ of a subspace of $\mathbb{R}^p$ with dimension $(\dim\cF^\prime\!-\!\dim\cF)$ such that $A^\prime\subset A$ and $\varphi^{-1}(A^\prime\times \{c\})$ transverse to $\cF$. Also, $\exists\ B^\prime \subset B$ such that $\psi^{-1}(B^\prime\times \{d\}))$ transverse to $\cF$.\smallskip \ \\
Let $g$ be the composition of the restriction map $f\vert_{A^\prime\times B^\prime\times C}$ and the quotient map as the diagram showed left below: 
$$
\xymatrix{
A\times B\times C\ar[r]^-{f}&f(A\times B\times C)\ar[d]^{\slash\sim}\\
A^\prime\times B^\prime \times C\ar@{^{(}->}[u]\ar[r]^-{g}&[f(A\times B\times C)].
}
\hspace{2cm}
\xymatrix{
\Hol(M,\cF^\prime)^{(1)}\ar[d]^{s}\ar[r]^-{\slash\sim}&\cH^{(1)}\ar[d]^{s}\\
M\ar[r]^{\pi}&W.
}
$$
We claim that $g$ is a local chart of $[\gamma]$, also the image of $g$ form a topological basis for the quotient topology on $\cH$. It suffices to prove the bijectivity of $g$. \\
\indent \underline{Injectivity}: If $g(a^\prime,b^\prime,c)=g(\tilde{a^\prime},\tilde{b^\prime},\tilde{c})$, that is $[f(a^\prime,b^\prime,c)]=[f(\tilde{a^\prime},\tilde{b^\prime},\tilde{c})]$, then, $s\big(f(a^\prime,b^\prime,c)\big)$ and $s\big(f(\tilde{a^\prime},\tilde{b^\prime},\tilde{c})\big)$ lie on the same leaf of $\cF$ and in the plaque $A\times C$. As $\varphi^{-1}(A^\prime\times C)$ is transverse to $\cF$, $a^\prime=\tilde{a^\prime}, c=\tilde{c}$. Also, $b^\prime=\tilde{b^\prime}$.\\
\indent \underline{Surjectivity}: For each $[\gamma]\in [f(A\times B\times C)]$, take $\gamma\!=\!f(x,y,z)$ a representative. $\exists !\ x^\prime\in A^\prime$ such that $\varphi^{-1}(x,z)$ and $\varphi^{-1}(x^\prime,z)$ lie on the same leaf of $\cF$. As the same, $\exists ! \ y^\prime\in B$ such that $r(f(x,y,z))$ and $r(f(x,y^\prime,z))$ lie on the same leaf of $\cF$. Clearly, $[\gamma]\!=\! [f(x^\prime,y^\prime,z)]\!=\![f(x,y,z)]$ with $(x^\prime,y^\prime,z)\!\in\! A^\prime\!\times\! B^\prime\!\times\! C$. \medskip
\ \\ 
\underline{Part II: $s,r$ are smooth submersions}
In the diagram right above, $s\!:\!\Hol(M,\cF^\prime)^{(1)}\!\rightarrow\! M$ and $\pi$ are both smooth submersions, so $s:\cH^{(1)}\rightarrow W$ is also a smooth submersion. The proof of $r$ is similar. \medskip
\ \\
\underline{Part III: Hausdorffness of $s$-fiber}\\
For $w\in W$, take $x\in \pi^{-1}(w)$. We know that the $s$-fiber $\Hol(M,\cF^\prime)_x$ is Hausdorff. To prove the Hausdorffness of $\cH_w$, it needs to prove the openness of the quotient map $\Hol(M,\cF^\prime)_x\rightarrow \cH_w$ and the closeness of the equivalence relation $\sim$ on $\Hol(M,\cF^\prime)_x\times \Hol(M,\cF^\prime)_x$. In \underline{Part I}, local charts $f$ of $\Hol(M,\cF^\prime)$ and $g$ of $\cH$ are open sets, so the quotient map is open. If $(\gamma_n,\delta_n)\rightarrow (\gamma,\delta)$ in $\Hol(M,\cF^\prime)_x\times \Hol(M,\cF^\prime)_x$ where $\forall n\in \mathbb{N}$, $\delta_n\sim\gamma_n$, then $\delta_n * \gamma_n^{-1}=1_{x}^{x}$. Passing to limit, $\delta * \gamma^{-1}=1_{x}^{x}$.\medskip
\ \\
\underline{Part IV The unit inclusion $i:W\rightarrow \cH^{(1)}$ is smooth}\\ 
Let $U$ be a small open of $W$ and $\sigma:U\rightarrow M$ be a smooth local section of $\pi:M\rightarrow W$. For $w\in U$,$i(w)=[\mathrm{1}_{\sigma(w)}]$ 
is smooth.\eb
\begin{Thm}\label{Thm1}
Let $\cF$ be the simple foliation of a fibration $M\rightarrow W$ of two compact manifolds. If $\cF^\prime$ is a foliation on $M$ with $\cF\subset \cF^\prime$, then the foliation $\cF^\prime$ projects to a foliation $\cF_W$ on $W$. Furthermore, if $\cF^\prime$ is a Riemannian foliation, $\cF_W$ also is.
\end{Thm}
\Preuve By proposition \ref{Prop4}, $\cH$ is a Lie groupoid. By \cite{Moe} Theorem 5.4 (iii), $\forall\ w\in W$, the orbit $\cH\!\cdot\! w$ is a submanifold of $W$. Followed by the construction of local chart of $\cH$, $\forall\ w\in W$, $\dim(\cH\!\cdot\! w)\!=\!(\dim\cF^\prime\!-\!\dim\cF)$. The distribution
$$
\bigcup_{w\in W}\big\{ T_w (\cH\!\cdot\! w)\vert w\in W \big\}
$$
is completely integrable and of constant rank, thus defines a foliation $\cF_W$ on $W$. As $\cH$ is the quotient of $\cF^\prime$, 
$T\cF^\prime\vert_m\rightarrow T\cF_W\vert_{\pi(m)}$ is surjective for all $m\in M$. The inheritance of Riemannian property will result from Proposition \ref{Prop14}. \eb
\begin{Rem}
The foliation $\cF_W$ can be regarded as the projection of $\cF^\prime$. The groupoid $\cH$ is exactly the holonomy groupoid of $\cF_W$.
\end{Rem}
\bigskip
\noindent The rest of this section is devoted to discuss (basic) differential forms. 
\begin{Def}
On a foliated manifold $(M,\cF)$, a differential form $\alpha\in \Omega^\bullet(M)$ of degree $\bullet$ is $\cF$-basic if $\mathcal{L}_X\alpha=i_X\alpha=0$ for all vector fields $X$ tangent to $\cF$, where $\mathcal{L}$ is Lie derivative and $i$ is contraction. Let $\big(\Omega^\bullet(M,\cF),d\big)$ denote the complex of basic forms on $(M,\cF)$ with the exterior differential $d$, $H^\bullet(M,\cF)$ denote its cohomology, called \textbf{basic cohomology}.
\end{Def}
\begin{Prop}\label{Prop13}
There is an identification of basic forms
$$
\Omega^\bullet(M,\cF^\prime)\simeq \Omega^\bullet(W,\cF_W).
$$
and a cohomological isomorphism
$$
H^\bullet(M,\cF^\prime)\simeq H^\bullet(W,\cF_W).
$$
\end{Prop}
\Preuve Let $\alpha\in \Omega^\bullet(M,\cF^\prime)$. Each $\cF^\prime$-basic form on $M$ is $\cF$-basic and identified to a form $\underline{\alpha}\in\Omega(W)$, \cite{Liu} Proposition 1.1.1. For each field $X$ on $W$ tangent to $\cF_W$, take a small open $U\subset W$. It suffices to calculate locally on $\pi^{-1}(U)$, for any lift $\widetilde{X}$ of $X$ tangent to $\cF^\prime$,  $$
i_X(\underline{\alpha})\simeq i_{\widetilde{X}}\alpha=0,\ 
\mathcal{L}_X(\underline{\alpha})\simeq \mathcal{L}_{\widetilde{X}}\alpha=0.
$$
The isomorphism of basic cohomology is clear.
\eb
Remark that a $\cF^\prime$-basic form on $M$ can be regarded as an element in $C^\infty\big(M,\Lambda(\nu\cF^\prime)^* \big)^{\cF^\prime\text{-}\mathrm{inv}}$, where $(-)^{\cF^\prime\text{-}\mathrm{inv}}$ means the $\cF^\prime$-invariance. 
\begin{Prop}\label{Prop14}
Generally, we have 
$$
C^\infty\big(M,S(\nu\cF^\prime)^*\oplus \Lambda(\nu\cF^\prime)^*  \big)^{\cF^\prime\text{-}\mathrm{inv}}\simeq C^\infty\big(W,S(\nu\cF_W)^*\oplus\Lambda (\nu\cF_W)^* \big)^{\cF_W\text{-}\mathrm{inv}},
$$
where $S$ is the symmetric algebra. In particular, 
$$
C^\infty\big(M,S^2(\nu\cF^\prime)^* \big)^{\cF^\prime\text{-}\mathrm{inv}}\simeq C^\infty\big(W,S^2(\nu\cF_W)^*\big)^{\cF_W\text{-}\mathrm{inv}}.
$$
That is, if $\cF^\prime$ is a Riemannian foliation, then $\cF_W$ also is. 
\end{Prop}
\section{Lift of Foliation}
We recall Molino's theory. Let $(M,\cF)$ be a Riemannian foliation. Its orthonormal frame bundle $O(\nu\cF)$ is naturally equipped with a lifted foliation $\tcF$. $\tcF$ is transversely parallelizable (T.P.) and its closure is a manifold $W$, called basic manifold. Remark that the bundle $O(\nu\cF)$ is quite specific, Molino constructed the lifted foliation $\tcF$ through a basic connection. But generally, the lift of foliation is preferred before construction of basic connection. This is the procedure in \cite{Liu}, where lifted foliations are generated by action of holonomy groupoid.\bigskip\ \\
Section 4.1 explains in detail a simple case, where $\cF$ partial transversely parallelizable (P.T.P) to $\cF^\prime$. The action of $\Hol(M,\cF^\prime)$ defines a lift $\tcF^\prime$ of $\cF^\prime$ in $O(\nu\cF)$. Furthermore, its projection on the basic manifold $W$ is a foliation. This is an application of Theorem \ref{Thm1}. Remark that the Riemannian assumption of $\cF^\prime$ is essential for getting rid of a pathological projection. A counterexample is in Remark \ref{Rem8}.\bigskip\ \\
In Section 4.2, we consider the non-P.T.P. case. One can first lift $\cF$ into the orthonormal frame bundle $O(\cF^\prime\slash\cF)$. Then, the lifted foliation $\hcF$ is partial T.P. to the pullback foliation $\hcF^\prime$. Remark that the Riemannian property of $\cF^\prime$ is not necessary for this step, but essential for a further work. Thus, the final form is Theorem \ref{Thm2}. 
\begin{Not}
Some notations in following text are listed below. For definition, see \cite{Mol} Chapter 2.
\begin{enumerate}
\item $\frak{X}(\cF)$ vector fields on $M$ tangent to $\cF$,$\frak{X}(M,\cF)$ Lie algebra of $\cF$-foliated vector fields on $M$;
\item $l(M,\cF)=\frak{X}(M,\cF)\slash\frak{X}(\cF)$ Lie algebra of transverse fields on $M$;
\item $\frak{X}(M,\cF,\cF^\prime)$ Lie algebra of vector fields both $\cF$ and $\cF^\prime$-foliated on $M$; 
\item $\frak{X}(\cF^\prime,\cF)$ Lie subalgebra of $\frak{X}(M,\cF,\cF^\prime)$ formed by $\cF$-foliated vector fields tangent to $\cF^\prime$; 
\item $l(\cF^\prime,\cF)\!=\!\frak{X}(\cF^\prime,\cF)\slash \frak{X}(\cF)$, $l(M,\cF,\cF^\prime)\!=\!\frak{X}(M,\cF,\cF^\prime)\slash \frak{X}(\cF)$;
\item $\xi\!=\!\cF^\prime\slash\cF$, $\perp:\!\nu\cF\rightarrow \nu\cF^\prime$ natural projection; 
\end{enumerate}
\end{Not}

\subsection{P.T.P. case} 
Let $\cF,\cF^\prime$ be Riemannian foliations on $M$, $\mathrm{Codim}\cF^\prime\!=\!q$, $\mathrm{Codim}\cF\!=\!q\!+\!s$. Take respectively transverse metrics $g_1,g_2$ of $\cF,\cF^\prime$, consider and fix $g=g_1\vee g_2$ as the transverse metric on $\nu\cF$. Let $\nu\cF\!\!=\!\xi\oplus \xi_\perp$ be the orthogonal decomposition with respect to $g$. A vector $v\in T\cF\vert_x$ with $x\in M$ is written by $v\!=\!v_{/\!/}\!+\!v_{\perp}$. An isometric isomorphism is $\sigma\!:\xi_\perp\rightarrow\!\nu\cF^\prime$. 
\begin{Def}
A foliation $\cF$ is \textbf{partial transversely parallelizable} (P.T.P.) to $\cF^\prime$ if there exist a family 
$$
Y_1,Y_2,\cdots,Y_s\in l(\cF^\prime,\cF),
$$
which forms a global frame for the bundle $\xi$.
\end{Def}
\begin{Def}
For $x,y\in M$, we define the linear map $\forall\ k\in \{1,2,\cdots,s\}$:
$$
|\!|_x^y:\xi\vert_x\rightarrow\xi\vert_y, Y_k\vert_x\mapsto Y_k\vert_y.
$$ 
\end{Def}
\noindent We may assume that $Y_k$ is orthonormal under $g$, thus $|\!|_x^y$ is an isometry. 
\begin{Def}\label{Def8}
A $\Hol(M,\cF^\prime)$-action on $\nu\cF$ is defined for each $\delta\in \Hol(M,\cF^\prime)$ and $v\in {\nu\cF}\vert_{s(\delta)}$, 
$$
\delta\!\cdot\! v= |\!|_{s(\delta)}^{r(\delta)} \, v_{/\!/}+\sigma^{-1}\big(\delta\!\cdot\!\sigma(v_\perp)\big).
$$
\end{Def}
\noindent It is easy to check that the $\Hol(M,\cF^\prime)$-action is isometric and respects the orthogonality $\xi\oplus\xi_\perp$. \\

\noindent Let $\cT$ be a complete transversal of $\cF$. By Corollary \ref{subtransversal}, there exist a complete subtransversal $\cU$ of $\cF^\prime$. An isomorphism is
\begin{equation}\label{ptpdecom}
(T\cT)\vert_\cU\simeq T\cU\oplus \xi\vert_\cU.
\end{equation}
Let $J^1(\cU)=\big\{(x,y,B)\mid x,y\in \cU,\ B:T_x\,\cU\rightarrow T_y\,\cU\ \mathrm{is\ }g_2\ \mathrm{isometric}\big\}$ be the isometric $1$-jet groupoid. There is an injective groupoid homomorphism
$$
j^\prime:\Hol(M,\cF^\prime)_\cU^\cU\rightarrow J^1(\cU).
$$
Thanks to Relation (\ref{ptpdecom}), the homomorphism
$$
J^1(\cU)\rightarrow J^1(\cT),\ B\mapsto B\oplus |\!|_{x}^{y}
$$
is well-defined. It implies an action of $J^1(\cU)$ on $(\nu\cF)\vert_\cU$ because $(\nu\cF)\vert_\cU\!\simeq\! (T\cT)\vert_\cU$. Furthermore, we have a composition of homomorphism 
\begin{equation}\label{eq02}
\Hol(M,\cF^\prime)_\cU^\cU\stackrel{j^\prime}{\longrightarrow} J^1(\cU)\rightarrow J^1(\cT)
\end{equation}
Let $\tM\!=\! O(\nu\cF)$ be the orthonormal transverse frame bundle of $\nu\cF$ and $\tcF$ be the lifted foliation by Molino's theory. Under the action of $\Hol(M,\cF^\prime)$ on $\nu\cF$, $\tM$ is $\Hol(M,\cF^\prime)$-equivariant by Remark \ref{Rem2.1}.

\begin{Prop}{\label{Prop6}}
The action of $\Hol(M,\cF^\prime)$ defines a lifted foliation $\tcF^\prime$ such that $(\tM,\tcF^\prime)$ is a principal $\cF$-bundle with $\tcF\subset \tcF^\prime$. Furthermore, $\tcF$ is P.T.P. to $\tcF^\prime$.
\end{Prop}
\Preuve First, $\tcF^\prime$ is well-defined. Followed by Relation \ref{eq02}, the closure of $\Hol(M,\cF^\prime)_\cU^\cU$ in $J^1(\cU)$ acts on $(T\cT)\vert_\cU\simeq \nu\cF\vert_\cU$. By Proposition \ref{Prop-2}, $(\tM,\tcF^\prime)$ is a principal $\cF$-bundle equipped with a $\tcF^\prime$-basic connection. The inclusion $\tcF\subset\tcF^\prime$ results from Proposition \ref{Prop1} (2). Since $Y_k\in l(\cF^\prime,\cF)$, the lifted field $\widetilde{Y_k}$ with respect to any $\tcF^\prime$-basic connection belongs to $l(\tcF^\prime,\tcF)$. Hence, $\tcF$ is P.T.P. to $\tcF^\prime$. \eb
\begin{Rem}
By Proposition \ref{Prop-1.1}, $\tcF$ and $\tcF^\prime$ are both Riemannian. Each $\tcF^\prime$-basic connection is automatically $\tcF$-basic. 
\end{Rem}
 
Let us recall the transverse parallelism defined in Section 3.3 \cite{Mol}. Let $\theta_T\in C^\infty\big(\tM,(\nu\tcF)^*\otimes \mathbb{R}^q\big)$ be the transverse fundamental form on $\tM$, \cite{Mol} Section 3.3. Recall the definition of $\theta_T$: for $z\in \tM$, $v\in \nu\tcF\vert_z$, 
$$
\theta_T(v)=z^{-1}\circ p_*(v),
$$
where $p_*:\nu\tcF\rightarrow \nu\cF$ is induced by the projection $p:\tM\rightarrow M$. Take $\omega$ a $\tcF^\prime$-basic connection.\\

\noindent A natural transverse parallelism $\big\{ \overline{\lambda_1},\cdots,\overline{\lambda_{q(q-1)/2}},\overline{u_1},\cdots,\overline{u_q}\big\}$ of $\nu\tcF$ is defined by the identities:
\begin{equation}\label{parallel}
\left\{
\begin{array}{lcl}
\omega(\overline{u_i})=0&;&\theta_T(\overline{u_i})=u_i\\
\omega(\overline{\lambda_j})=\lambda_j&;&\theta_T(\overline{\lambda_j})=0\\
\end{array}
\right.
\end{equation}
where $\{\lambda_1,\lambda_2,\cdots,\lambda_{q(q-1)/2}\}$ is a basis of $so(q)$ and $\{u_1,\cdots, u_q\}$ is the canonical basis of $\mathbb{R}^q$.
\begin{Lemme}\label{Lemme5}
Let $\widetilde{Y_k}$ be the lift of $Y_k$ with respect to $\omega$, for $k\!=\!1,\cdots, s$. The form $\theta_T$ is $\tcF^\prime$-invariant. In particular, $\mathcal{L}_{\widetilde{Y_k}} \theta_T=0$. 
\end{Lemme}
\Preuve It suffices to show the invariance of $\theta_T$ under $\Hol(\tM,\tcF^\prime)$-action. By Proposition \ref{Prop-1.5} and the groupoid isomorphism    
$$
\Hol(M,\cF^\prime)\ltimes \tM\simeq \Hol(\tM,\tcF^\prime), 
$$
for $(\gamma,z)\in \Hol(M,\cF^\prime)\ltimes \tM$, $v\in \nu\tcF\vert_z$, the action $(\gamma,z)\cdot v$ is realized by the direct sum of the horizontal lift and the identity map on vertical subbundle. Thus, we have an identity:    
$$
p_*\big( (\gamma,z)\cdot v\big)=\gamma\cdot (p_* v).
$$
The target of $(\gamma,z)$ is the isometry $w\in \tM$: 
$$
w:\mathbb{R}^q\stackrel{z}{\rightarrow} \nu\cF\vert_{p(z)}\stackrel{\gamma\cdot}{\rightarrow} \nu\cF\vert_{p(w)}
$$
A brief calculation proves the lemma. 
\[
\begin{array}{rl}
\big((\gamma,z)^*\theta_T\big)(v)&=\theta_T\vert_w\big((\gamma,z)\cdot v\big)=w^{-1}\circ p_*\big((\gamma,z)\cdot v\big)\\ \\[-0.5em]
&=w^{-1}\circ \gamma\cdot (p_*(v))=z^{-1}\circ p_*(v)=\theta_T(v).
\end{array}
\]
\eb
\begin{Prop}
The transverse parallelism of $\nu\tcF$ in Relation (\ref{parallel}) belongs to $l(\tM,\tcF,\tcF^\prime)$. \end{Prop}
\Preuve Since the actions of $\Hol(M,\cF^\prime)$ and $SO(q)$ on $\tM$ commute, $\tcF$ and $\tcF^\prime$ are both $SO(q)$-invariant. Thus, $\overline{\lambda_j}\in l(\tM,\tcF,\tcF^\prime)$. Let us prove $\overline{u_i}\in l(\tM,\tcF,\tcF^\prime)$.\\
It suffices to show $[\widetilde{Y_k},\overline{u_i}]\!=\!0$ in $l(\tM,\tcF)$ because $\widetilde{Y_k}$ is a parallelism of $\tcF^\prime\slash\tcF$. By Lemma \ref{Lemme5}, 
$$
\theta_T\big([\widetilde{Y_k},\overline{u_i}]\big)=\big(\mathcal{L}_{\widetilde{Y_k}}\theta_T\big)(\overline{u_i})-\widetilde{Y_k}\big(\theta_T(\overline{u_i})\big)=-\widetilde{Y_k}(u_i)=0,
$$
$$
\omega\big([\widetilde{Y_k},\overline{u_i}]\big)=\widetilde{Y_k}\big(\omega(\overline{u_i})\big)-\big(\mathcal{L}_{\widetilde{Y_k}}\omega\big)(\overline{u_i})=0,
$$
which implies $[\widetilde{Y_k},\overline{u_i}]=0$ in $l(\tM,\tcF)$.\eb
\begin{Cor}\label{Cor3}
The parallelism $\widetilde{Y_k}$ of $\tcF^\prime\slash \tcF$ can be completed to a transverse parallelism of $\nu\tcF$ by $\widetilde{Y_l}\in l(\tM,\tcF,\tcF^\prime)$. Furthermore, the quotient fields $(\widetilde{Y_l})_\perp\!\in\! l(\tM,\tcF^\prime)$ form a transverse parallelism of $\nu\tcF^\prime$, $\tcF^\prime$ is transversely parallelizable.
\end{Cor}
\subsection{Lifted foliations and commuting sheaf }
\indent Let $\ta$ be the commuting sheaf of $(\tM,\tcF)$, \cite{GT}. This section is devoted to discuss the relation between $\tcF^\prime$ and $\ta$. We use the notation presented in \cite{Mol} Section 4.4. Let $U$ be an connected open subset of $\tM$, denote $\mathcal{C}(U)$ the set of local transverse fields which commute with $l(\tM,\tcF)$. It is a subalgebra of $l\big(U,\tcF\vert_U\big)$. Next proposition is critical.
\begin{Prop}\label{Lemme7}
If $Y\in \mathcal{C}(U)$ with $Y\vert_{m_0}\in (\tcF^\prime\slash\tcF)\vert_{m_0}$ for some $m_0\in U$, then $\forall \ m\in U$, $Y\vert_m\in (\tcF^\prime\slash\tcF)\vert_m$. Consequently, the distribution 
$$
\displaystyle\bigcup_{m\in\tM}\Big((\ta\!\cdot\!\tcF)\vert_m\cap \tcF^\prime\vert_m\Big)
$$
is of constant rank. 
\end{Prop}
\Preuve The proof is analogous to the proof of  the Hausdorffness of $\mathcal{C}(M,\cF)$ in \cite{Mol} Page 126. It suffices to prove the openness and closeness in $U$ of the set
$$
\Big\{m\in U,\ Y\vert_m\in (\tcF^\prime\slash\tcF)\vert_m \Big\}.
$$
The closeness is evident. Since $Y$ commutes with global transverse fields, in particular, for $k\!=\!1,\!\cdots\! ,s$, $[Y,\widetilde{Y_k}]=0$, then $Y$ preserves the foliation $\tcF^\prime$. Hence, its projection $Y_\perp\in l\big(U,\tcF^\prime\vert_U\big)$.\\
Take $S$ a transversal of $\tcF^\prime$ at $m_0$. On $S$, under identification, $Y_\perp$ commutes with $(\widetilde{Y_l})_\perp$ in Corollary \ref{Cor3}. Hence, $Y_\perp\vert_{m_0}=0$ implies that $Y_\perp$ is zero in an open neighborhood of $m_0$. \eb 
\begin{Rem}\label{Rem8}
If $\tcF^\prime$ is non-Riemannian, the rank $(\ta\!\cdot\!\tcF)\vert_m\cap \tcF^\prime\vert_m$ may vary along $m$. We have an example inspired by Example \ref{Ex1}: let $M\!=\!S^1\!\times\! S^1\!\times\! S^1\!\times\! S^1$ with coordinates $(x,y,z,t)$. Define foliations with $\alpha\in\mathbb{R}\backslash \mathbb{Q}$ by \\[-1em]
\[
\begin{array}{rl}
\cF\!&\!=\!\big<\partial\slash\partial x\!+\!\alpha\,\partial\slash\partial t\big>\\ \\[-1em]
\cF^\prime\!&\!=\!\big<\partial\slash\partial x\!+\!\alpha\,\partial\slash\partial t,\partial\slash\partial x\!+\!f(y)\,\partial\slash\partial z\big>,
\end{array}
\]
where $f$ is a non-constant function on $S^1$. Evidently, $\cF\subset \cF^\prime$. Since $\cF$ is the Kronecker foliation (Killing), 
$\fa$ is exactly the center of $l(M,\cF)$. Take $\partial\slash\partial t$ as a representative generator of $\fa$. Denote $m=(x,y,z,t)$, we see
\begin{itemize}
\item[$\bullet$] if $y\in S^1$ with $f(y)=0$, $\partial\slash\partial t\vert_m\in \cF^\prime\vert_m$
\item[$\bullet$] if $y\in S^1$ with $f(y)\neq 0$, $\partial\slash\partial t\vert_m\notin\cF^\prime\vert_m$.
\end{itemize}
Therefore, the rank of $(\ta\!\cdot\!\tcF)\vert_m\cap \tcF^\prime\vert_m$ varies with $m$. On the other hand, if $\tcF^\prime$ is Riemannian, $\cF^\prime$ is affirmatively Riemannian. Consequently, the Riemannian assumption of $\cF^\prime$ is indispensable.
\end{Rem}
\begin{Prop}\label{Prop11}
The distribution 
$$
\ta\!\cdot\!\tcF^\prime=\displaystyle\bigcup_{x\in\tM} \Big((\ta\!\cdot\!\tcF)\vert_x + \tcF^\prime\vert_x\Big)
$$
defines a Riemannian foliation on $\tM$.
\end{Prop}
\Preuve \underline{Involutiveness} The rank  
$$
\mathrm{Rank}\big(\ta\!\cdot\!\tcF^\prime\big)\vert_x=\mathrm{Rank}\big(\tcF^\prime\vert_x\big)+\mathrm{Rank}(\ta\vert_x)-\mathrm{Rank}\big((\ta\!\cdot\!\tcF)\vert_x\cap \tcF^\prime\vert_x\big)
$$
is constant. Check locally the involutiveness of $\ta\!\cdot\!\tcF^\prime$. Equivalently, we may assume that $\ta$ is globally constant. For $Y_1,Y_2\in \ta$, $Z_1,Z_2\in l(\tcF^\prime,\tcF),$ 
$$
\big[Y_1\!+\!Z_1,Y_2\!+\!Z_2\big]\!=\!\big[Y_1,Y_2\big]\!+\!\big[Y_1,Z_2\big]\!+\!\big[Z_1,Y_2\big]\!+\!\big[Z_1,Z_2\big]\!=\!\big[Z_1,Z_2\big]\in l(\tcF^\prime,\tcF).
$$   
\noindent \underline{Riemannian property} Consider the transverse metric $\tilde{g}^\prime$ on $\nu\tcF^\prime$ by disposing the parallelism $(\widetilde{Y_l})_\perp$ orthonormal. Locally, $\forall\ Y\in \ta$,
{\footnotesize
$$
(\mathcal{L}_Y \tilde{g}^\prime)\big((\widetilde{Y_l})_\perp,(\widetilde{Y_{l^\prime}})_\perp\big)\!=\!Y\Big(\tilde{g}^\prime\big((\widetilde{Y_l})_\perp,(\widetilde{Y_{l^\prime}})_\perp\big)\Big)-\tilde{g}^\prime\big([Y,\widetilde{Y_l}]_\perp,(\widetilde{Y_{l^\prime}})_\perp\big)-\tilde{g}^\prime\big((\widetilde{Y_l})_\perp,[Y,\widetilde{Y_{l^\prime}}]_\perp\big)\!=\!Y\Big(\tilde{g}^\prime\big((\widetilde{Y_l})_\perp,(\widetilde{Y_{l^\prime}})_\perp\big)\Big).
$$
}
$\tilde{g}^\prime\big((\widetilde{Y_l})_\perp,(\widetilde{Y_{l^\prime}})_\perp\big)$ is a $\tcF^\prime$-basic function, so $\tcF$-basic. Each $\tcF$-basic function is $\overline{\tcF}$-basic, that is $\ta\!\cdot\!\tcF$-basic. So, $Y\big(\big(\tilde{g}^\prime(\widetilde{Y_l})_\perp,(\widetilde{Y_{l^\prime}})_\perp\big)\big)=0$ Hence, $\tilde{g}^\prime$ is $\ta$-invariant.\\ \\
Consider the projection $\tau:\nu\tcF^\prime\rightarrow \nu(\ta\!\cdot\!\tcF^\prime)$. By Proposition \ref{Lemme7}, its kernel $\mathrm{Ker}(\tau\vert_x)$ is of constant rank, thus $\mathrm{Ker}(\tau)$ is a subbundle of $\nu\tcF^\prime$.  Hence, $\tilde{g}^\prime$ defines an orthogonal decomposition
$$
\nu\tcF^\prime=\mathrm{Ker}(\tau)\oplus \mathcal{N}.
$$
Each element in $\ta$ locally preserve $\tcF^\prime$ and $[\ta,\ta]\subset \ta$, the subbundle $\mathrm{Ker}(\tau)$ is stable under $\Hol(\tM,\tcF^\prime)$-action as well as $\ta$-action. The $\ta$ and $\Hol(\tM,\tcF^\prime)$ invariance of $\tilde{g}^\prime$ inherits to $\mathcal{N}$. Followed by the identification $\mathcal{N}\simeq \nu(\ta\!\cdot\!\tcF^\prime)$, we obtain a $\ta$ and $\Hol(\tM,\tcF^\prime)$-invariant metric on $\nu(\ta\!\cdot\!\tcF^\prime)$. In other words, $\ta\!\cdot\!\tcF^\prime$ is a Riemannian foliation. \eb
\subsection{Non-P.T.P. Case}
Let the manifold $M$, the foliations $\cF,\cF^\prime$, the metric $g$, the complete transversal (resp. subtransversal) $\cT$ (resp. $\cU$) be same as Section 5.1. In this part, $\cF$ is not supposed to be P.T.P. to $\cF^\prime$. \medskip \ \\
Let $\hM=O(\xi)$ be the orthonormal frame bundle of $\xi$, $\hat{p}$ denote the projection $\widehat{M}\rightarrow M$. By Proposition \ref{Prop1}(1), \ref{Prop1.1}, Remark \ref{Rem2.1}, the action of $\Hol(M,\cF)$ defines a lifted foliation $\hcF$ on $\hM$ such that $(\hM,\hcF)$ is a principal $\cF$-bundle. Denote by $\hcT=\hM\vert_\cT$. Next proposition results from Propositions \ref{Prop-2}, \ref{Prop-1}, \ref{Prop0}.
\begin{Prop}\label{Prop0.1}
The principal bundle $(\hM,\hcF)$ is a $\cF$-bundle. There is a groupoid isomorphism
$$
\overline{\Hol(M,\cF)_\cT^\cT}\ltimes \hcT\simeq \overline{\Hol(\hM,\hcF)_\hcT^\hcT}.
$$
In particular, there is an identification between Molino sheaf of $(M,\cF)$ and $(\hM,\hcF)$.
\end{Prop}
\begin{Def}
We define $\widehat{\cF^\prime}=\hat{p}^{-1}(\cF^\prime)$ to be the \textbf{pullback foliation} on $\widehat{M}$, endowed with the pullback transverse metric $\hat{p}^{*}g_{2}$. 
\end{Def}
\begin{Prop}
The foliation $\widehat{\cF}$ is P.T.P. to $\widehat{\cF^\prime}$.
\end{Prop}
\Preuve Clearly, $\widehat{\cF}\subset \widehat{\cF^\prime}$. Define analogously the fundamental $\xi$-transverse form $\theta^\prime_T\in \mathcal{A}^1(\widehat{M},\mathbb{R}^s)$ by 
$$
{\theta^\prime_T}\vert_{z}(v)=z^{-1}\big(\overline{ T_z\hat{p}(v)  } \big)_{/\!/}, \ v\in T_z\tM
$$
where $/\!/:\nu\cF\rightarrow \xi$ is the orthogonal projection and in this proof $\overline{(\ )}$ is the projection $TM\rightarrow \nu\cF$. Consider $\{\lambda_1,\lambda_2,\cdots,\lambda_{s(s-1)/2}\}$ the canonical basis of $\frak{so}(s)$ and $\{u_1,\cdots, u_s\}$ the canonical basis of $\mathbb{R}^s$. Then the $\hcF$-basic fields 
$$
\{\overline{\lambda_1},\overline{\lambda_2},\cdots,\overline{\lambda_{s(s-1)/2}},\overline{u_1},\cdots, \overline{u_s}\}
$$
defined by the identities:
\[
\left\{
\begin{array}{lcl}
\widehat{\omega}(\overline{u_i})=0&;&\theta^\prime_T(\overline{u_i})=u\\
\widehat{\omega}(\overline{\lambda_j})=\lambda&;&\theta^\prime_T(\overline{\lambda_j})=0\\
\overline{u_i}\in l(\cF^\prime,\cF)
\end{array}
\right.
\]
form a transverse parallelism of $\widehat{\cF}_2/\widehat{\cF}_1$. Therefore, $\widehat{\cF}$ is P.T.P. to $\widehat{\cF^\prime}$.\eb

\noindent  Apply Proposition \ref{Prop6} to $(\widehat{M}, \widehat{\cF}_1,\widehat{\cF}_2)$. Let $\tM\!=\!O(\nu\widehat{\cF})$ be the orthonormal transverse frame of $(\widehat{M},\widehat{\cF})$ with respect to $\widehat{g}$ and $W\!=\!\tM\slash\overline{\tcF}$ be the basic manifold.
\begin{Prop}
We have the lift of foliations $\tcF, \tcF^\prime$ in $\tM$ with $\tcF\subset \tcF^\prime$ such that $(\tM,\tcF^\prime)$ is a principal $\cF$-bundle equipped with a $\tcF^\prime$-basic connection. Furthermore, $\tcF$ and $\tcF^\prime$ are both transversely parallelizable and $\tcF$ is P.T.P. to $\tcF^\prime$. 
\end{Prop}
\begin{Prop}
The foliation $\tcF^\prime$ projects to a Riemannian foliation $\cF_W$ on $W$.
\end{Prop}
\Preuve The projection of vectors tangent to $\tcF^\prime$ by $\pi$ is as same as vectors tangent to the foliation $(\ta\!\cdot\!\tcF^\prime)$. By Proposition \ref{Prop11}, $(\ta\!\cdot\!\tcF^\prime)$ is a Riemannian foliation. Followed by Theorem \ref{Thm1}, the proof is complete.\eb
\noindent In summary, we have the following theorem.
\begin{Thm}\label{Thm2}
Let $M$ be a foliated manifold equipped with Riemannian foliations
$\cF, \cF^\prime$ satisfying $\cF\subset \cF^\prime$. Let $\widehat{M}\!=\!O(\cF^\prime\slash\cF)$ be the orthonormal frame bundle of $\cF^\prime\slash\cF$ and $\hcF^\prime$ be the pullback foliation on $\hM$.
\begin{enumerate}
\item $\widehat{M}$ is a principal $\cF$-bundle with a lifted Riemannian foliation $\hcF$. The foliation $\hcF$ is partial transversely parallelizable to $\hcF^\prime$;
\item $\tM\!=\!O(\nu\widehat{\cF})$ is a principal $\cF$-bundle with lifted Riemannian foliations $\tcF,\tcF^\prime$. The foliations $\tcF,\tcF^\prime$ are transversely parallelizable in $\tM$ and $\tcF$ is partial transversely parallelizable to $\tcF^\prime$;
\item $\tcF^\prime$ projects to a Riemannian foliation $\cF_W$ on $W$. 
\end{enumerate}
See the diagram below.
\begin{equation}\label{Diagram1}
\xymatrix{
&(\tM,\tcF,\tcF^\prime)\ar[r]^-{\pi}\ar[d]^-{\widetilde{p}}&(W,\cF_W)\\
(M,\cF,\cF^\prime)&(\hM,\hcF,\hcF^\prime)\ar[l]^-{\widehat{p}}&
}
\end{equation}  
\end{Thm}
To finish this section, we propose the identification of Molino sheaves. Let $\fa,\ha$ and $\ta$ denote respectively Molino sheaf on $(M,\cF)$, $(\widehat{M},\widehat{\cF})$ and $(\tM,\tcF)$. Apply Proposition \ref{Prop0}, we obtain
\begin{Prop}\label{Prop14}
There is a Molino sheaf identification
$$
\fa\simeq \ha \simeq \ta.
$$
\end{Prop}
\noindent In following text, we do not distinguish $\fa,\ha$ and $\ta$.
\section{Cohomological Isomorphisms}
\noindent Let us recall the definition of transverse action of a Lie algebra.
\begin{Def}
A \textbf{transverse action} of a Lie algebra $\frak{h}$ on a foliated manifold $(M,\cF)$ is a Lie homomorphism $\frak{h}\rightarrow l(M,\cF),\ Y\mapsto Y_M$. If $Y_M\!=\!0$ at some point of $M$ implies $Y\!=\!0$, we say the action is free.
\end{Def}
In this section, assume that in Theorem \ref{Thm2}, $\cF$ is a {\bf Killing foliation}. A Killing foliation is a foliation whose Molino sheaf is globally constant, that is, $\ta$ is exactly the center of transverse fields, 
\begin{equation}\label{eq09}
\ta=\mathrm{Center}\big(l(\tM,\tcF)\big).
\end{equation}
It is an Abelian Lie algebra, called \textbf{Molino algebra}. Followed by Relation (\ref{eq09}), there is a transverse action
\begin{equation}\label{eq10}
\fa\rightarrow l(\tM,\tcF)^{SO(\hat{q})}
\end{equation}
where $\hat{q}=\mathrm{Codim}\hcF$, $(-)^{SO(\hat{q})}$ means $SO(\hat{q})$-invariance. 
By Proposition \ref{Lemme7}, if $Y\in \fa$ such that $Y$ belongs to $\tcF^\prime\slash\tcF$ at one point, $Y$ belongs to $l(\tcF^\prime,\tcF)$. Then, the transverse metric of $\tcF$ defines an orthogonal decomposition of Lie algebra 
$$
\fa=\fa_{\perp}\oplus \fa_{/\!/}.
$$
Naturally, the transverse action (\ref{eq10}) is decomposed into two free actions
$$
\fa_{\perp}\rightarrow l(\tM,\tcF^\prime)^{SO(\hat{q})}\ \ \  \fa_{/\!/}\rightarrow l(\tcF^\prime,\tcF)^{SO(\hat{q})}.
$$
Concentrate on the first one, by equivariance, it induces transverse actions
$$
\fa_{\perp}\rightarrow l(\hM,\hcF^\prime)^{SO(s)}\ \ \fa_{\perp}\rightarrow l(M,\cF^\prime),
$$
where $s\!=\!\dim\cF^\prime-\dim\cF$. For the notion of equivariant basic cohomology with $C^\infty$-coefficients, see \cite{DKV}. Next proposition is a generalization of \cite{GT} Proposition 4.9, which is the case $\cF^\prime=\cF$.
\begin{Prop}
We have the following cohomological isomorphism
\begin{equation}\label{eq11}
H^\infty_{\fa_{\perp}}(M,\cF^\prime)\simeq H^\infty_{\fa_{\perp}\times so(\hat{q})}(\tM,\tcF^\prime)\simeq H^\infty_{so(\hat{q})}(W,\cF_W).
\end{equation}
\end{Prop}
\Preuve 
Since on $\Omega(\tM,\tcF^\prime)$, the action of $\fa_{\perp}$ and $SO(\hat{q})$ are commutative and free,
$$
H^\infty_{\fa_{\perp}}(\hM,\hcF^\prime)\simeq H^\infty_{\fa_{\perp}\times so(\hat{q})}(\tM,\tcF^\prime)\simeq H^\infty_{so(\hat{q})}(\tM,\fa_\perp\!\!\cdot\!\!\tcF^\prime).
$$
The foliation $\fa_\perp\!\cdot\!\tcF^\prime=\fa\!\cdot\!\tcF^\prime$. Followed by Proposition \ref{Prop13},
$$
\Omega^\bullet(\tM,\fa_\perp\!\cdot\!\tcF^\prime)=\Omega^\bullet(\tM,\fa\!\cdot\!\tcF^\prime)\simeq \Omega^\bullet(W,\cF_W).
$$
As $(\tM,\fa\!\cdot\!\tcF^\prime)\rightarrow (W,\cF_W)$ is $SO(\hat{q})$-equivariant, we have
$$
H^\infty_{so(\hat{q})}(\tM,\fa_\perp\!\cdot\!\tcF^\prime)\simeq H^\infty_{so(\hat{q})} (W,\cF_W).
$$
Followed from that $\hcF^\prime$ is the pullback foliation of $\cF^\prime$,
$$
H^\infty_{\fa_\perp}(\hM,\hcF^\prime)\simeq H^\infty_{\fa_\perp}(M,\cF^\prime).  
$$
\eb  
\section{Geometric realization of $\hat{A}$-genus character}
This section aims to give a geometric realization of cohomological isomorphism (\ref{eq11}) through equivariant basic $\hat{A}$-genus characters. See \cite{GS} for the equivariant de Rham theory, \cite{BGV} for the notion of equivariant $\hat{A}$-genus character. Review some results in \cite{Liu1}. \medskip \ \\
Let $E$ be a vector bundle associated to a principal $\cF$-bundle $(P,\cF_P)$ of structure group $G$, $\frak{h}$ be a Lie algebra. If there exist a diagram of Lie homomorphism
$$
\xymatrix{
&l(P,\cF_P)^G\ar[d]\\
\frak{h}\ar[r]\ar[ru]&l(M,\cF)
}
$$
and a $\frak{h}$-invariant basic connection $\omega_P$ on $P$, then by the associated $\cF$-basic connection $\nabla^E$ and the action of $\frak{h}$ on $E$, the $\frak{h}$-equivariant curvature $R^E(Y)=R^E\!+\!\mu^E(Y)$ is well-defined where $\mu$ denotes moment map. The form given by $\mathrm{Tr}(e^{R^E(Y)})$ defines the equivariant basic Chern character. In this work, we study the form
$$
\hat{A}_{\frak{h}}(\nabla^E)=\det{}^{\!\frac{1}{2}}\left[\dfrac{R^{E}(Y)/2}{\sinh(R^{E}(Y)/2)}\right]
$$
which defines the equivariant basic $\hat{A}$-genus character. As the normal bundle of a Riemannian foliation is always associated to its transverse orthonormal frame bundle, which is a $\cF$-bundle, an associated basic connection always exists. The invariance of connection on normal bundle under a transverse action is equivalent to that on principal bundle. Therefore, in following texts, it suffices to check the invariance of connection on normal bundle without passing to principal bundle.\medskip \ \\
\textbf{Correspondence under $H_{so(\hat{q})}^{\infty}(W,\cF_W) \simeq H_{\fa_\perp\times so(\hat{q})}^{\infty}(\tM,\tcF^\prime)$}\smallskip \
\begin{Lemme}\label{Lemme12}
Let $G$ be a Lie group. Let $E, F$ be two $G$-vector bundles (vector bundle equipped with an action of $G$) on $M$. Let $\Phi:E\rightarrow F$ be a $G$-equivariant bundle isomorphism, that is, $\forall\ g\in G$, $g\circ \Phi=\Phi\circ g$. If $\nabla^F$ is a $G$-invariant connection on $F$ and $\nabla^E=\Phi^*\nabla^F$, the $G$-equivariant curvatures $R^F(X)$, $R^E(X)$ satisfy
$$
R^E(X)=\Phi^*\big(R^F(X)\big).
$$ 
\end{Lemme}
\Preuve Clearly, $R^E\!=\!\Phi^*R^F$. The moment maps satisfy $\mu^E(X)\!=\!\Phi^*\big(\mu^F(X)\big)$ followed by the $G$-equivariance of $\Phi$. By definition, $R^F(X)\!\!=\!\!R^F\!+\!\mu^F(X)$.  \eb

Take a $SO(\hat{q})$-invariant basic connection $\nabla^{\nu\cF_W}$ on $\nu\cF_W$ (it always exists because $SO(\hat{q})$ is compact), then $SO(\hat{q})$-equivariant basic form
$$
\hat{A}_{so(\hat{q})}(\nabla^{\nu\cF_W})=\det{}^{\!\frac{1}{2}}\left[\dfrac{R^{\nu\cF_W}(X)/2}{\sinh(R^{\nu\cF_W}(X)/2)}\right],
$$
is well-defined where $R^{\nu\cF_W}(X)$ is the $SO(\hat{q})$-equivariant curvature. Its cohomological class $\hat{A}_{so(\hat{q})}(\nu\cF_W)$, independent from the choice of $\nabla^E$, belongs to $H_{so(\hat{q})}^\infty(W,\cF_W)$.\\
Let $\mathcal{V}$ denote the trivial bundle
generated by $\fa_\perp$. Since $\nu\tcF^\prime$ is trivializable, there exist a subbundle $\cU$ such that
$$
\nu\tcF^\prime=\mathcal{U}\oplus \mathcal{V},
$$
where $\mathcal{U}$ (resp. $\mathcal{V}$) possesses a trivialization $\{Y_l\in l(\tM,\tcF^\prime)\}$ (resp. $\{Y_k\in l(\tM,\tcF^\prime)\}$) by Corollary \ref{Cor3}.\\
Clearly, $Y_l$ projects to $\underline{Y_l}\in l(W,\cF_W)$, $(\underline{\ })$ denotes the projection. Thus, $\nu\cF_W$ is also trivializable. Consider the pullback bundle $\pi^*(\nu\cF)$ where $\pi^*(\underline{Y_l})$ can be identified to $Y_l$. 
\begin{Lemme}\label{Lemme13}
A $SO(\hat{q})$-equivariant isomorphism $\Phi:\mathcal{U}\rightarrow \pi^*(\nu\cF_W)$ is given by $Y_l\simeq \pi^*(\underline{Y_l})$.
\end{Lemme}
\Preuve The action of $SO(\hat{q})$ on $\pi^*(\nu\cF_W)$ comes from $\nu\cF_W$. $\forall\ g\in SO(\hat{q})$, the relation $\underline{g\cdot Y_l}=g\cdot \underline{Y_l}$ results from the $SO(\hat{q})$-equivariance of $\pi:(\tM,\tcF^\prime)\rightarrow (W,\cF_W)$. \eb 
Take a $SO(\hat{q})$-invariant $\cF_W$-basic connection $\nabla^{\nu\cF_W}$ on $\nu\cF_W$, let $\nabla^{\pi^*\nu\cF_W}=\pi^*\nabla^{\nu\cF_W}$ be its pullback which is naturally $(\fa_\perp\!\!\cdot\!\! \tcF^\prime)$-basic. By Lemma \ref{Lemme13}, $\nabla^\cU=\Phi^*\nabla^{\pi^*\nu\cF_W}$ is $\fa_\perp\!\!\times\!\! SO(\hat{q})$-invariant. It is easy to see that the moment map of $\fa_\perp$ is zero.\\
Let $\nabla^{\mathcal{V}}$ be the trivial connection defined by $\nabla^{\mathcal{V}}_Z(fY_k)=Z(f)Y_k$. For $X\in so(\hat{q})$, $X_{\tM}$ is the vector field on $\tM$ generating by $X$. Lie derivative is $\mathcal{L}^\mathcal{V}_X Y_k=[X_{\tM},Y_k]=0$. It is easy to check that $\nabla^{\mathcal{V}}$ is $(\fa_\perp\!\!\cdot\!\!\tcF^\prime)$-basic and $\fa_\perp\!\!\times\!\! SO(\hat{q})$-invariant. The moment map of $\fa_\perp\!\!\times\!\! SO(\hat{q})$ on $\mathcal{V}$ is zero. 
\begin{Prop}
The inclusion $i_{so(\hat{q})}$ induces the following equivariant basic $\hat{A}$-genus character relation
\[
\begin{array}{rcl}
i_{so(\hat{q})}:H_{so(\hat{q})}^{\infty}(W,\cF_W)&\rightarrow& H_{\fa_\perp\times so(\hat{q})}^{\infty}(\tM,\tcF^\prime)\\
\hat{A}_{so(\hat{q})}(\nu\cF_W)&\mapsto&\hat{A}_{\fa_\perp\times so(\hat{q})}(\nu\tcF^\prime)
\end{array}
\]
\end{Prop}
\Preuve  Consider the $(\fa_\perp\!\!\cdot\!\!\tcF^\prime)$-basic $\fa_\perp\!\!\times\!\! SO(\hat{q})$-invariant connection $\nabla^{\nu\tcF^\prime}\!=\!\nabla^\cU\oplus\nabla^\cV$.
Then, the $\fa_\perp\!\!\times\!\! so(\hat{q})$-equivariant basic curvature $R^{\nu\tcF^\prime}(Y,X)\!=\!\big(R^\cU\!+\!\mu^{\cU}(X)\big)\oplus R^{\cV}\!=\!R^\cU(X)\oplus R^\cV$. Notice that $R^\cV\!=\!0$, we calculate explicitly
\[
\begin{array}{rl}
\hat{A}_{\fa_\perp\!\times\! so(\hat{q})}\big(\nabla^{\nu\tcF^\prime}\big)=&
\det{}^{\!\frac{1}{2}}\left[\dfrac{R^\cU(X)/2}{\sinh(R^\cU(X)/2)}\right]\wedge \det{}^{\!\frac{1}{2}}\left[\dfrac{R^\cV/2}{\sinh(R^\cV/2)}\right]\\ \\[-0.5em]
=&\pi^*\det{}^{\!\frac{1}{2}}\left[\dfrac{R^{\nu\cF_W}(X)/2}{\sinh(R^{\nu\cF_W}(X)/2)}\right] \\ \\[-0.5em]
=&\pi^* \big(\hat{A}_{so(\hat{q})}(\nabla^{\nu\cF_W})\big)=i_{so(\hat{q})}\big(\hat{A}_{so(\hat{q})}(\nabla^{\nu\cF_W})\big).
\end{array}
\]
\eb \medskip\ \\
\textbf{Correspondence under $H_{\fa_\perp}^{\infty}(\hM,\hcF^\prime)\simeq H_{\fa_\perp\times so(\hat{q})}^{\infty}(\tM,\tcF^\prime)$}\medskip \ \\
The action of $\fa_\perp$ on $\nu\hcF^\prime$ (resp. $\nu\tcF^\prime$) is Lie bracket because the action of $\fa$ on $\nu\hcF$ (resp. $\nu\tcF$) is. Without ambiguity, in this part, let $\cU$ (resp. $\cV$) denote the horizontal lift of $\nu\hcF^\prime$ in $\nu\tcF^\prime$ with respect to a fixed $\tcF^\prime$-basic connection (resp. the canonical vertical subbundle). Followed by relation (\ref{parallel}), there exist a parallelism $\{\overline{u_i}\}$ of $\cU$ and $\overline{\lambda_j}$ of $\cV$. 
\begin{Def}
We define a bundle isomorphism
\[
\begin{array}{rcl}
\Psi:\cU&\rightarrow&p^*(\nu\hcF^\prime)\\
\overline{u_i}&\mapsto&\widetilde{u_i}
\end{array}
\]
where $\widetilde{u_i}\vert_{\widetilde{m}}\!=\!p_*\big(\overline{u_i}\vert_{\tilde{m}}\big)$ with $\widetilde{m}\in \tM$ and $p_*$ is vector projection.
\end{Def}
\begin{Prop}\label{Prop25}
The isomorphism $\Psi$ is $\fa_\perp\!\times\! so(\hat{q})$-equivariant.
\end{Prop}
\Preuve $\bullet$ $\fa_\perp$-equivariance.  For $Y\in\fa_\perp$, any representative $\widehat{Y}\in \frak{X}(\hM,\hcF^\prime)$ with lift $\widetilde{Y}\in \frak{X}(\tM,\tcF^\prime)$,
\[
\begin{array}{rl}
\big(\mathcal{L}^{p^*\nu\hcF^\prime}_{\widetilde{Y}} \widetilde{u_i}\big)\vert_{\widetilde{m}}=&\frac{d}{dt}_{\vert_{t=0}} (e^{-t\widetilde{Y}})_*\big(\widetilde{u_i}\vert_{e^{t\widetilde{Y}}\cdot\widetilde{m}}\big)=\frac{d}{dt}_{\vert_{t=0}} (e^{-t\widetilde{Y}})_*\circ p_*\big(\overline{u_i}\vert_{e^{t\widetilde{Y}}\cdot\widetilde{m}}\big)\\ \\
=&p_*\Big(\frac{d}{dt}_{\vert_{t=0}}(e^{-t\widehat{Y}})_* \big(\overline{u_i}\vert_{e^{t\widetilde{Y}}\cdot \widetilde{m}}\big) \Big)=p_*\Big(\big(\mathcal{L}^\cU_{\widetilde{Y}}\overline{u_i} \big)\vert_{\widetilde{m}}\Big)  
\end{array}
\]
$\bullet$ $SO(\hat{q})$-equivariance.  The action of $SO(\hat{q})$ on $p^*(\nu\hcF^\prime)$ and on $\cU$ is respectively
$$
(g\!\cdot\! \overline{u_i})\vert_{\widetilde{m}}=g_*\big(\overline{u_i}\vert_{\widetilde{m}g^{-1}}\big), \ (g\!\cdot\! \widetilde{u_i})\vert_{\widetilde{m}}=\widetilde{u_i}\vert_{\widetilde{m}g^{-1}}=p_*(\overline{u_i}\vert_{\tilde{m}g^{-1}}).
$$
Thus, $\Psi\big((g\!\cdot\! \overline{u_i})\vert_{\widetilde{m}}\big)=p_*\!\circ\! g_*\big(\overline{u_i}\vert_{\widetilde{m}g^{-1}}\big)=p_*\big(\overline{u_i}\vert_{\widetilde{m}g^{-1}}\big)=(g\!\cdot\! \widetilde{u_i})\vert_{\widetilde{m}}$.\eb
It is easy to see that $\Psi$ preserve $\hcF$-holonomy invariance. Followed by Proposition \ref{Prop1.1}, $\nu\hcF^\prime$ satisfies \textbf{Assumption} for the action of $\Hol(\hM,\hcF)$, then equipped with a $\fa_\perp$-invariant basic connection $\nabla^{\nu\hcF^\prime}$. Its pullback $\nabla^{p^*\nu\hcF^\prime}\!\!=\!\!p^*\nabla^{\nu\hcF^\prime}$ is naturally $SO(\hat{q})$-invariant. Besides, the $\fa_\perp$-invariance of $\nabla^{p^*\nu\hcF^\prime}$ results from the $\fa_\perp$-equivariance of $p:(\tM,\tcF^\prime)\rightarrow (\hM,\hcF^\prime)$. By Proposition \ref{Prop25}, the connection 
$$
\nabla^\cU=\Psi^*\big(\nabla^{p^*\nu\hcF^\prime}\big)
$$
is $\tcF^\prime$-basic and $\fa_\perp\!\!\times\!\! SO(\hat{q})$-invariant. Therefore, the moment map $\mu^\cU(Y)=\mathcal{L}^\cU_{\widetilde{Y}}-\nabla^\cU_{\widetilde{Y}}$ (resp. $\mu^{p^*\nu\hcF^\prime}$, $\mu^{\nu\hcF^\prime}$) is well-defined, independent of representative $\widetilde{Y}$.\\
For $X\in so(\hat{q})$, $\mu^{p^*\nu\hcF^\prime}(X)=0$, so $\mu^\cU(X)\!=\!\Psi^*\mu^{p^*\nu\hcF^\prime}(X)\!=\!0$. Next lemma is clear.
\begin{Lemme}
For $X\in so(\hat{q})$, $Y\in\fa_\perp$, we have $R^\cU(X,Y)=\hat{p}^*R^{\nu\hcF^\prime}(Y)$.
\end{Lemme}
Let $\nabla^\cV$ be the trivial connection on $\cV$ defined by $\nabla^\cV_Z(f\overline{\lambda_j})=Z(f)\overline{\lambda_j}$. 
For $X\in so(\hat{q})$, the moment map $\mu^\cV(X)(\overline{\lambda_j})=\overline{(\mathrm{ad}X)(\lambda_j)}$ because $\nabla^\cV_{X_\tM}\overline{\lambda_j}=0$ and the relation $\mathcal{L}^\cV_X \overline{\lambda_j}=\overline{(\mathrm{ad}X)(\lambda_j)}$ comes from 
\[
\left\{\begin{array}{l}
\omega\big(\mathcal{L}^\cV_X \overline{\lambda_j}\big)=(\mathrm{ad}X)(\lambda_j)\\ \\[-1em]
\theta_T\big(\mathcal{L}^\cV_X \overline{\lambda_j}\big)=0.
\end{array}\right.
\]
The matrix of $\mu^\cV(X)$ is similar to $\mathrm{ad}(X)$.
\begin{Prop}
The connection $\nabla^\cV$ is $\hcF^\prime$-basic and $\fa_\perp\!\!\times\!\! SO(\hat{q})$-invariant. 
\end{Prop}
\Preuve The basic fields $\overline{\lambda_j}$ are $\tcF^\prime$-holonomy invariant, so the connection $\nabla^\cV$ is $\hcF^\prime$-basic.\\
$$
\mathcal{L}^\cV_{\widetilde{Y}}\big(\nabla^\cV_Z(f\overline{\lambda_j})\big)=\widetilde{Y}\big(Z(f)\big)\overline{\lambda_j}+Z(f)[\widetilde{Y},\overline{\lambda_j}]=\widetilde{Y}\big(Z(f)\big)\overline{\lambda_j}. 
$$   
$$
\nabla^\cV_{[\widetilde{Y},Z]}(f\overline{\lambda_j})+\nabla^V_Z\big([\widetilde{Y},f\overline{\lambda_j}]\big)=[\widetilde{Y},Z](f)\overline{\lambda_j}+Z\big(\widetilde{Y}(f)\big)\overline{\lambda_j}=\widetilde{Y}\big(Z(f)\big)\overline{\lambda_j}.
$$      
The $\fa_\perp$-invariance is proved, the $SO(\hat{q})$-invariance is similar.\eb
\noindent The connection $\nabla^{\nu\tcF^\prime}\!=\!\nabla^{\cU}\!\oplus\!\nabla^{\cV}$ fits to define the following equivariant basic form
$$
\hat{A}_{\fa_\perp\times so(\hat{q})}(\nabla^{\nu\tcF^\prime})(Y,X)=\det{}^{\!\frac{1}{2}}\left[\dfrac{R^{\nu\tcF^\prime}(Y,X)/2}{\sinh(R^{\nu\tcF^\prime}(Y,X)/2)}\right].
$$  
\noindent Let $\hat{A}_{so(\hat{q})}$ denote the function $\det{}^{\!\frac{1}{2}}\left[\dfrac{(\mathrm{ad}X)/2}{\sinh((\mathrm{ad}X)/2)}\right]$ of $so(\hat{q})$.
\begin{Prop}
The inclusion $i_{\fa_\perp}:H_{\fa_\perp}^{\infty}(\hM,\hcF^\prime)\rightarrow H_{\fa_\perp\times so(\hat{q})}^{\infty}(\tM,\tcF^\prime)$ induces the following equivariant basic $\hat{A}$-genus character relation
$$
i_{\fa_\perp}\big(\hat{A}_{\fa_\perp}(\nu\cF^\prime)\big)=\hat{A}^{-1}_{so(q)}\cdot\hat{A}_{\fa_\perp\times so(\hat{q})}(\nu\tcF^\prime).
$$
\end{Prop}
\Preuve We calculate explicitly 
\[
\begin{array}{rl}
R^{\nu\tcF^\prime}(Y,X)=&\big[(R^{\cU}\!+\!\mu^\cU(Y)\!+\!\mu^\cU(X)\big]\oplus \big[R^\cV\!+\!\mu^\cV(Y)\!+\!\mu^\cV(X)\big]\\ \\[-0.8em]
=&\big[R^\cU\!+\!\mu^\cU(Y)\big]\oplus \mu^\cV(X).
\end{array}
\]
$$
\begin{array}{rl}
\hat{A}_{\fa_\perp\times so(q)}\big(\nabla^{\nu\tcF^\prime}\big)=&\det{}^{\!\frac{1}{2}}\left[\dfrac{R^\cU(Y)/2}{\sinh(R^\cU(Y)/2)}\right]\wedge \det{}^{\!\frac{1}{2}}\left[\dfrac{\mu^\cV(X)/2}{\sinh(\mu^\cV(X)/2)}\right]\\ \\[-0.8em]
=&\det{}^{\!\frac{1}{2}}\left[\dfrac{R^{p^*\nu\hcF^\prime}(Y)/2}{\sinh(R^{p^*\nu\hcF^\prime}(Y)/2)}\right]\cdot \det{}^{\!\frac{1}{2}}\left[\dfrac{(\mathrm{ad}X)/2}{\sinh((\mathrm{ad}X)/2)}\right]\\ \\[-0.8em]
=&\hat{A}_{so(\hat{q})}\cdot p^*\det{}^{\!\frac{1}{2}}\left[\dfrac{R^{\nu\hcF^\prime}(Y)/2}{\sinh(R^{\nu\hcF^\prime}(Y)/2)}\right]=\hat{A}_{so(\hat{q})}\cdot p^*\hat{A}_{\fa_\perp}(\nabla^{\nu\hcF^\prime}). \\ \\[-0.8em]
\end{array}
$$
\eb
Finally, as $\hcF^\prime$ is pullback foliation, there is a natural isomorphism $p^*(\nu\cF^\prime)\simeq \nu\hcF^\prime$. 
\begin{Prop}
The pullback $p^*: H^\infty_{\fa_\perp}(M,\cF^\prime)\rightarrow H^\infty_{\fa_\perp}(\hM,\hcF^\prime)$ induces an equivariant basic $\hat{A}$-genus character relation
$$
p^*\big(\hat{A}_{\fa_\perp}(\nu\cF^\prime)\big)=\hat{A}_{\fa_\perp}(\nu\hcF^\prime)
$$
\end{Prop}
\Preuve The $\fa_\perp$-actions on $\nu\cF^\prime$ and $\nu\hcF^\prime$ are both Lie bracket, so compatible. Take a $\cF^\prime$-basic $\fa_\perp$-invariant connection $\nabla^{\nu\cF^\prime}$ and consider $p^*\nabla^{\nu\cF^\prime}$.
$$
\hat{A}_{\fa_\perp}(p^*\nabla^{\nu\cF^\prime})=p^*\hat{A}_{\fa_\perp}(\nabla^{\nu\cF^\prime}).
$$
\eb
The following theorem is a summary.
\begin{Thm}\label{Thm3}
Let $M$ be a foliated manifold equipped with a Killing foliation $\cF$ with Molino algebra $\fa$ and a Riemannian foliation $\cF^\prime$ satisfying $\cF\subset \cF^\prime$. Let $\fa_\perp$ be the Lie subalgebra orthogonal to $\cF^\prime$ of $\fa$. Then, there is a cohomological isomorphism as follows. Furthermore, the equivariant basic $\hat{A}$-genus characters of normal bundle give a geometric realization in following way
\[
\begin{array}{ccc}
H^\infty_{\fa_\perp}(M,\cF^\prime)&\simeq&H^\infty_{so(\hat{q})}(W,\cF_W)\\ \\[-0.5em]
\hat{A}_{\fa_\perp}(\nu\cF^\prime)&\simeq&\hat{A}^{-1}_{so(\hat{q})}\cdot\hat{A}_{so(q)}(\nu\cF_W),
\end{array}
\]
where \small $\hat{A}_{so(\hat{q})}=\det{}^{\!\frac{1}{2}}\left[\dfrac{(\mathrm{ad}X)/2}{\sinh((\mathrm{ad}X)/2)}\right]$\normalsize  is the function of $so(\hat{q})$ with a constant $\hat{q}$.
\end{Thm}
\begin{Cor}
In Molino's theory, for a Killing foliation $(M,\cF)$ of codimension $q$, its lift foliation $(\tM,\tcF)$ with Molino algebra $\fa$ and the basic manifold $W$, under the equivariant cohomological isomorphism, we have the following equivariant (basic) $\hat{A}$-genus character relation
\[
\begin{array}{rcl}
H^\infty_{\fa}(M,\cF)&\simeq&H^\infty_{so(\hat{q})}(W)\\ \\[-0.5em]
\hat{A}_{\fa}(\nu\cF)&\simeq&\hat{A}^{-1}_{so(q)}\cdot\hat{A}_{so(q)}(W),
\end{array}
\]
where $\hat{A}_{so(q)}(W)$ is the $so(q)$-equivariant $\hat{A}$-genus of the tangent bundle $TW$.
\end{Cor}

\end{document}